\documentclass[reqno,12pt]{amsart}

\usepackage{amssymb} \usepackage{amsfonts} \usepackage{amsmath}
\usepackage{amsthm} \usepackage{epsfig} 
\usepackage{color}
\usepackage{amscd}
\usepackage[all]{xy}
\usepackage{booktabs}
\usepackage{arydshln}

\newtheorem{lemma}{Lemma}[section]
\newtheorem{teo}[lemma]{Theorem}
\newtheorem{prop}[lemma]{Proposition}
\newtheorem{cor}[lemma]{Corollary}

\theoremstyle{definition}

\theoremstyle{remark}
\newtheorem{rem}[lemma]{Remark}

\newcommand{\matQ} {\ensuremath {\mathbb{Q}}}
\newcommand{\matZ} {\ensuremath {\mathbb{Z}}}

\newcommand{\matP} {\ensuremath {\mathbb{P}}}

\newcommand{\matRP} {\ensuremath {\mathbb{RP}}}

\newcommand{\nota} [1] {\caption{\footnotesize{#1}}}
\newcommand{\matr} [4] {\left(\begin{array}{@{}c@{\ }c@{}} #1 & #2 \\ #3 & #4 \\ \end{array} \right)}

\newfont{\Got}{eufm10 scaled 1200}

\newcommand{\GL}{{\rm GL}}

\newcommand{\cref}{c^{\rm ref}}

\newcommand{\Iso}{{\rm Isom}}

\newcommand{\timtil}{\begin{picture}(12,12)
\put(2,0){$\times$}\put(2,4.5){$\sim$}\end{picture}}

\newcommand{\lens}[2]{L({ #1},{ #2})}
\newcommand{\seifeul}[2]{\big(#1,#2\big)}
\newcommand{\seifuno}[3]{\big(#1,({#2},{ #3})\big)}
\newcommand{\seifunoeul}[4]{\big(#1,({ #2},{ #3}),#4\big)}
\newcommand{\seifdue}[5]{\big(#1,({#2},{ #3}),
                       ({#4},{#5})\big)}
\newcommand{\seifdueeul}[6]{\big(#1,({ #2},{ #3}),
                       ({#4},{ #5}),#6\big)}
\newcommand{\seiftre}[7]{\big(#1,({ #2},{ #3}),
                       ({ #4},{ #5}),
                       ({ #6},{ #7})\big)}

\newcommand{\bigu}[4]{\bigcup\nolimits_{{\tiny{\matr {#1} {#2} {#3} {#4}}}\phantom{\Big|}\!\!}}
\newcommand{\bigb}[4]{\big/_{{\tiny{\matr {#1} {#2} {#3} {#4}}}\phantom{\Big|}\!\!}}

\newcommand{\alzamolto}{\phantom{\Bigg(}\!\!\!\!\!}
\newcommand{\alza}{\phantom{\Big(}\!\!\!\!\!}
\newcommand{\alzapoco}{\phantom{(}\!\!\!\!\!}

\author{Bruno Martelli}
\address{Dipartimento di Matematica ``Tonelli'', Largo Pontecorvo 5, 56127 Pisa, Italy}
\email{martelli at dm dot unipi dot it}

\title{Dehn surgery \\ on the minimally twisted seven-chain link}

\begin{document}

\begin{abstract}
We classify all the exceptional Dehn surgeries on the minimally twisted chain links with six and seven components.
\end{abstract}

\maketitle

\section*{Introduction}
Let $M$ be a compact orientable three-manifold with some boundary tori. We say as usual that $M$ is \emph{hyperbolic} if its interior admits a finite-volume complete hyperbolic metric (which is then unique by Mostow and Prasad's rigidity theorem). Recall that a \emph{Dehn filling} of $M$ is the operation that consists of attaching solid tori to some (possibly all) of the boundary components of $M$, a manipulation that is essentially determined by the choice of some \emph{slopes} in the chosen boundary tori. 

\begin{figure}
\begin{center}
\includegraphics[width = 12 cm]{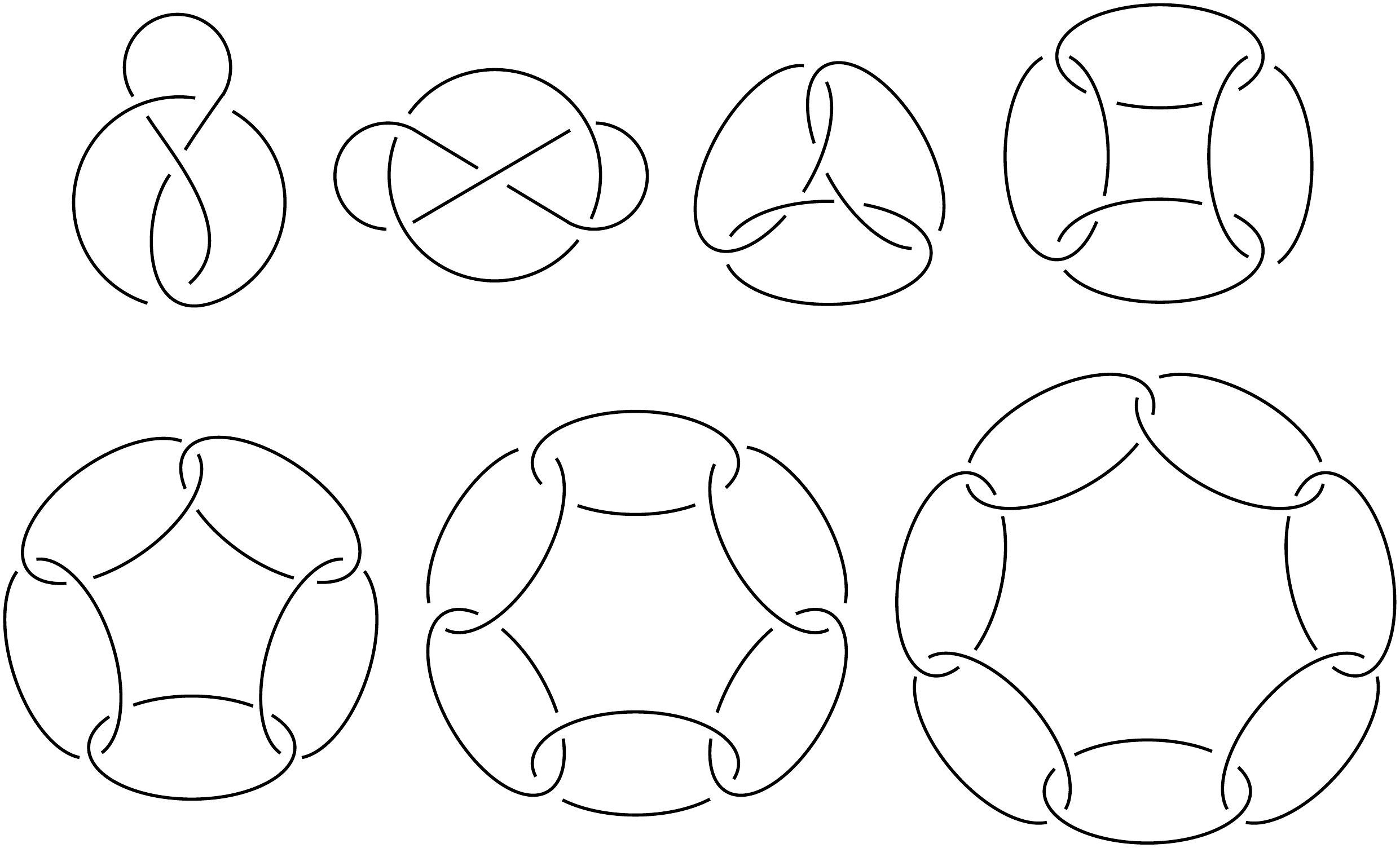}
\nota{A notable sequence of hyperbolic links with $i\leq 7$ components. These are the \emph{figure-eight knot}, the \emph{Whitehead link}, and some chain links with $3,\ldots, 7$ components. Those with $5, 6,$ and $7$ components are \emph{minimally twisted}.}
\label{new_sequence:fig}
\end{center}
\end{figure}

We say that a Dehn filling is \emph{hyperbolic} if the resulting manifold is still hyperbolic, and \emph{exceptional} otherwise. The goal of this paper is to make a further step in the classification of all the exceptional fillings in a natural sequence of hyperbolic link complements, initiated in \cite{MaPe} and \cite{MaPeRo}. The sequence is shown in Figure \ref{new_sequence:fig}. The main result is Theorem \ref{main:teo} below, where we exhibit a complete classification of all the exceptional fillings of the last two link complements shown in the figure.

\subsection*{The sequence}
\begin{figure}
\begin{center}
\includegraphics[width = 7 cm]{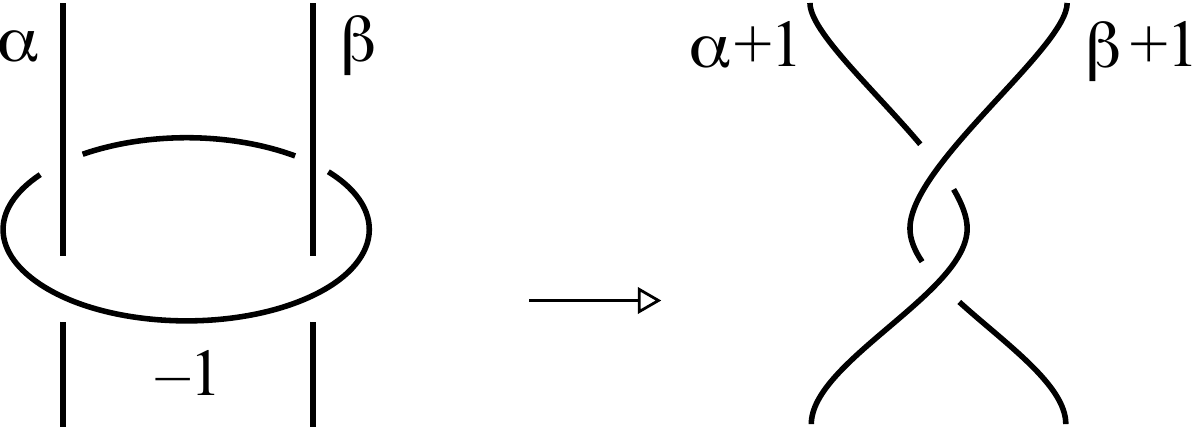}
\nota{A blow-down.}
\label{Rolfsen:fig}
\end{center}
\end{figure}

Figure \ref{new_sequence:fig} contains a sequence of notable hyperbolic links. These are the \emph{figure-eight knot}, the \emph{Whitehead link}, and some particular chain links with $i=3,\ldots, 7$ components. Let $M_1,\ldots, M_7$ be the complements of these links. Each $M_i$ is conjectured \cite{Ago:min} to have smallest volume among hyperbolic manifolds with $i$ cusps (this conjecture has been proved in the cases $i=1, 2$, and $4$ by Cao -- Meyerhoff \cite{CaMe}, Agol \cite{Ago:min}, and Yoshida \cite{Yo}). Another important feature of this sequence is that each $M_i$ is obtained as a $(-1)$-filling of the subsequent one $M_{i+1}$, as one sees via a blow-down as in Figure \ref{Rolfsen:fig}. See \cite{Thu, KPR, NR} for more information on hyperbolic chain links.

The manifolds $M_1,\ldots, M_7$ appear naturally in many contexts. The manifold $M_3$ was called the \emph{magic manifold} by Gordon and Wu in \cite{GoWu} because of its many interesting fillings; it plays a role in the study and/or classification of the closed hyperbolic manifolds of smallest volume
\cite{GMM}, of the pseudo-Anosov mapping classes with small dilatation \cite{Hi, KKT}, and of the cusped hyperbolic 3-manifold with the largest number of exceptional fillings \cite{LM}. 

\begin{figure}
\begin{center}
\includegraphics[width = 12 cm]{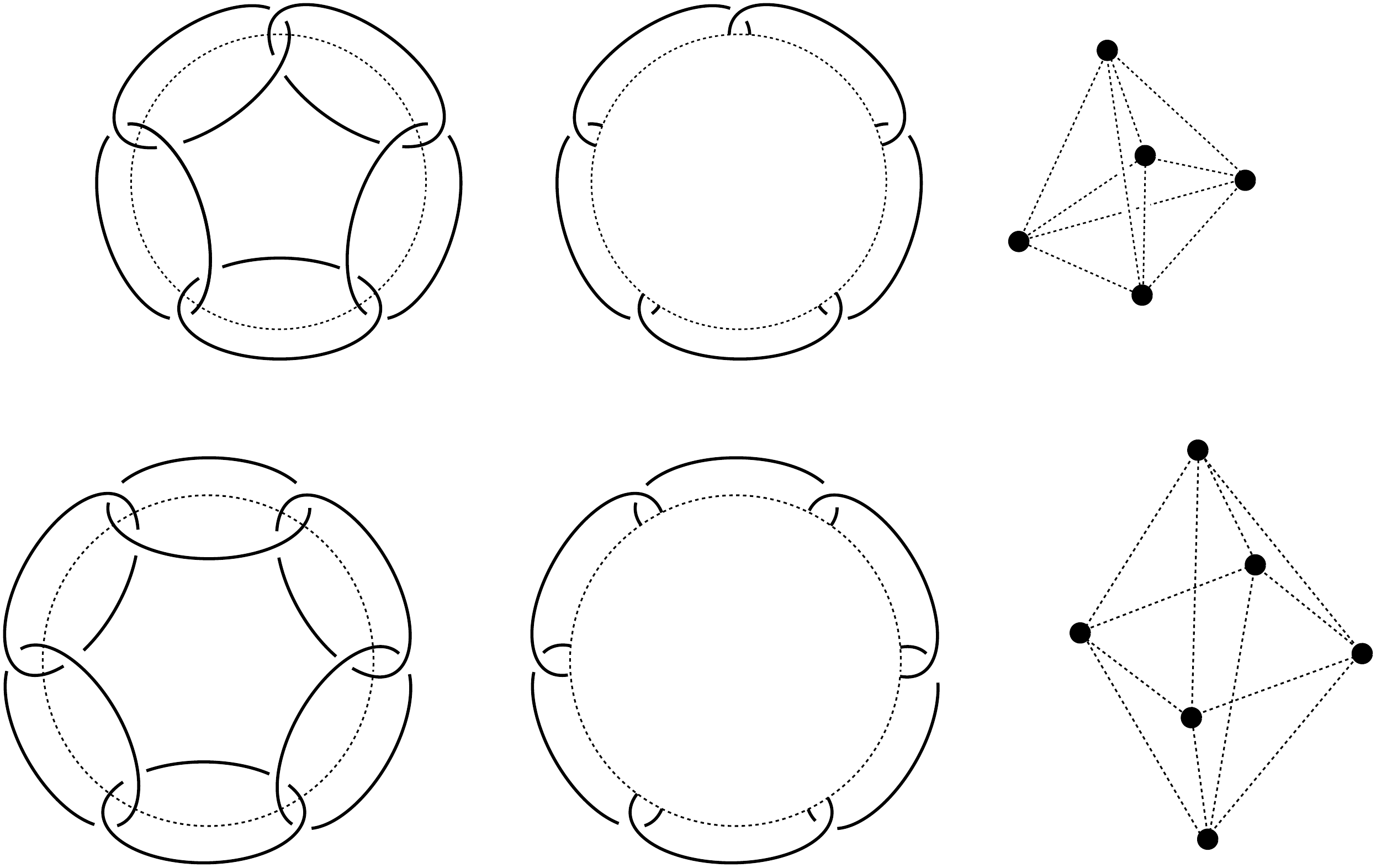}
\nota{The manifolds $M_5$ and $M_6$ double cover two natural orbifolds. The deck transformation is a $\pi$ rotation around the dotted axis.}
\label{graphs:fig}
\end{center}
\end{figure}

The manifolds $M_4, M_5$, and $M_6$ form a superb triple of highly symmetric hyperbolic manifolds. They decompose into regular ideal octahedra, tetrahedra, and octahedra respectively, and they are characterised (together with $M_3$) by the configuration of the many thrice-punctured spheres they contain \cite{Yo2}. The manifolds $M_5$ and $M_6$ cover two very natural hyperbolic orbifolds, shown in Figure \ref{graphs:fig}. The first is the boundary of the 5-simplex and decomposes into 5 regular ideal tetrahedra. The second is obtained by mirroring an ideal regular octahedron. The manifolds $M_3$, $M_5$, and $M_6$ are principal congruence link complements \cite{BR}, while $M_5$ and $M_6$ are also the smallest hyperbolic 3-manifolds admitting a regular tessellation \cite{G}.

The fillings of $M_4$ were used to classify the four-manifolds with shadow-complexity one \cite{KMN} and to build knots with long unknotting tunnels \cite{CFP}.
It was noted in \cite{Thu2, DuTh} that many cusped manifolds in the census \cite{CHW} are obtained by filling $M_5$. Among these, we find many Berge knots complements \cite{Bak} and other hyperbolic manifolds with interesting exceptional surgeries \cite{ALR, BDH, EM, GoLu, Kang}. 

The question of classifying all the exceptional fillings of $M_6$ was raised in \cite{EM2}. We answer to this question here. The manifold $M_6$ appears in the construction of hyperbolic four-manifolds with arbitrarily many cusps \cite{KM} and of arithmetic link extensions \cite{Bak2} in $S^3$. Its fundamental group is biorderable \cite{KR}. By filling $M_6$ we obtain yet more hyperbolic manifolds with interesting exceptional fillings \cite{EM2, Te}.

The manifold $M_7$ lacks all the beautiful symmetries of $M_4, M_5$, and $M_6$. It is the first and the only non-arithmetic manifold in the sequence $M_1, \ldots, M_7$.

\subsection*{The exceptional Dehn fillings}
The hyperbolic Dehn fillings of the figure-eight knot complement $M_1$ were famously described by Thurston in his notes \cite{Thu}. The exceptional fillings of the magic manifold $M_3$ were then classified by Martelli and Petronio in \cite{MaPe}. 
Later on, all the exceptional fillings on $M_5$ were listed by Martelli, Petronio, and Roukema in \cite{MaPeRo}. 

The main result of this paper is a complete classification of all the exceptional filings of the complement $M_7$ of the minimally twisted chain link with seven components, the last one of the sequence in Figure \ref{new_sequence:fig}. This of course includes also a classification of all the exceptional fillings of $M_6$.

\subsection*{Acknowledgements} The author would like to thank Nathan Dunfield and the anonymous referee for providing helpful comments on a first draft of this paper.

\section{Main results}

\subsection{The general strategy}
The exceptional Dehn fillings on a multi-cusped hyperbolic manifold $M$ may be infinite in number, but can typically be described using a finite amount of data. For instance, the magic manifold contains infinitely many exceptional fillings, which could be grouped into explicit families and were fully described in few pages in \cite{MaPe}.

To classify the exceptional fillings of $M_7$ we adopt the same general strategy of \cite{MaPe, MaPeRo}, that we now outline. More generally, let $M$ be any hyperbolic manifold with some boundary tori. Recall that a filling of $M$ is determined by a set of slopes, one for each filled boundary torus (we are allowed to leave some boundary tori unfilled). Following \cite{MaPeRo}, we say that an exceptional Dehn filling on $M$ is \emph{isolated} if any proper subset of the chosen slopes produces a hyperbolic Dehn filling. Thurston's Dehn filling Theorem implies the following:

\begin{teo} 
Every hyperbolic $M$ has only finitely many isolated exceptional fillings.
\end{teo}
To classify all the exceptional fillings of $M$ we must fulfill the following tasks:
\begin{enumerate}
\item classify all the isolated exceptional fillings of $M$,
\item recognise the filled manifolds in each case, and 
\item if necessary, proceed recursively on each filled manifold.
\end{enumerate}
The third point is necessary if the filled exceptional manifold contains some hyperbolic piece in its prime or JSJ decomposition, a case that did not occur in \cite{MaPe} and \cite{MaPeRo}, but that will arise here in this paper.

We could achieve all these objectives for the complement $M_7$ of the minimally twisted chain link with seven components. To this purpose we made an essential use of the formidable programs SnapPy \cite{Sna}, Regina \cite{Reg}, and Recognizer \cite{Rec}. Task (1) was fulfilled via {\tt find\_exceptional\_fillings.py}, a python script written by the author already used in \cite{MaPeRo} and publicly available \cite{M} to be performed on any cusped hyperbolic three-manifold. The computer-assisted proof is rigorous thanks to the {\tt hikmot} libraries \cite{hikmot}.

\subsection{The output}
When accomplished, the general strategy produces finitely many families of exceptional fillings, but as the number of cusps increases their number explodes, and it soon becomes impossible to write fully comprehensive tables as it was done for the magic manifold in \cite{MaPe}. 

We now describe the outcome of our research. The most concise amount of useful information that we can give is the following. 

\setlength{\tabcolsep}{10 pt}

\begin{table}
\begin{center}
\begin{tabular}{lrrrrrrrr}
\hline
& \multicolumn{5}{c}{Number of cusps filled} & \\ \cline{2-8}
Manifold & 1 & 2 & 3 & 4 & 5 & 6 & 7 & Total\\ \hline
$M_1$ & 10 &     &     &    & & & & 10 \\
$M_2$ & 12 & 14  &    &     & & & & 26 \\
$M_3$ & 15 & 15 & 52 & &  & & & 82 \\
$M_4$ & 16 & 24 & 96 & 492 & & & & 628 \\
$M_5$ & 15 & 30 & 180 & 780 & 4818 & & & 5823 \\
$M_6$ & 12 & 30 & 240 & 1572 & 7080 & 46680 & & 55614 \\
$M_7$ & 14 & 14 & 91 & 987 & 7119 & 32977 & 214007 & 255209 \\
\hline
\alzapoco
\end{tabular}
\end{center}
\nota{Numbers of isolated exceptional fillings on $M_i$.}
\label{exceptional_1:table}
\end{table}

\begin{teo} \label{main:teo}
The number of isolated exceptional fillings of $M_1,\ldots, M_7$ is shown in Table \ref{exceptional_1:table}. The complete lists of fillings can be downloaded from \cite{M}.
\end{teo}

The first impression that we get from looking at Table \ref{exceptional_1:table} is that the number of isolated exceptional fillings of $M_i$ with a fixed number $k$ of slopes is roughly constant as $i=1,\ldots, 7$ varies, and grows roughly exponentially in $k$.  

\subsection{Data reduction}
The manifold $M_7$ has 255,209 isolated exceptional fillings overall, and we would like to describe what these filled manifolds are. We cannot of course describe them all on a single table; instead, we try to reduce the amount of data that is necessary to understand and present them in some reasonable way.

It was already remarked in \cite{MaPeRo} that all the exceptional fillings of $M_5$ can actually be deduced from a very short lists of rules: a collection of 7 basic exceptional fillings plus a list of 5 isometries of $M_5$ and of some of its fillings generate all the exceptional fillings. Everything could be described in \cite{MaPeRo} in a half-page long theorem. (See Remark \ref{correction:rem} below for some corrections of the tables in \cite{MaPeRo}.)

We would like to find a similar small generating set of rules for the manifolds $M_6$ and $M_7$. As a first step, we quotient the exceptional fillings of $M_i$ by the action of its isometry group $\Iso(M_i)$. The isometry groups of $M_1,\ldots, M_7$ have order:
$$8, \qquad 8,\qquad 12,\qquad 64,\qquad 240,\qquad 192,\qquad 28.$$
These are respectively
$$D_8, \qquad D_8, \qquad D_{12}, \qquad G_{64}, \qquad S_5 \times \matZ_{2}, \qquad D_8 \times S_4, \qquad D_{28}.$$
Here $D_{2n}$ is the dihedral group of order $2n$ and the symbol $G_{64}$ indicates some non-abelian group of order $64$ that does not split as a direct product.
The manifolds $M_4, M_5$, and $M_6$ are arithmetic, decompose into regular ideal tetrahedra or octahedra, and have an extraordinary number of symmetries. On the other hand, the last manifold $M_7$ is not arithmetic \cite{NR} and has only few isometries: the 28 symmetries of the chain link that one infers from Figure \ref{new_sequence:fig} and nothing more than that. 

\begin{table}
\begin{center}
\begin{tabular}{lrrrrrrrr}
\hline
& \multicolumn{5}{c}{Number of cusps filled} & \\ \cline{2-8}
Manifold & 1 & 2 & 3 & 4 & 5 & 6 & 7 & Total\\ \hline
$M_1$ & 6    &     &     &    & & & & 6\\
$M_2$ & 6    & 8  &     &    & & & & 14\\
$M_3$ & 5    & 3  & 14 &  & & & & 22\\
$M_4$ & 2    & 2  &   4   &  22 & & & & 30\\
$M_5$ & 1    & 1  & 3  & 7 & 48 & & & 60 \\
$M_6$ & 1   & 2  & 4 & 22 & 79 & 529 & & 637 \\
$M_7$ & 2  & 2 & 9 & 73 & 522 & 2362 & 15357 & 18327 \\
\hline
\alzapoco
\end{tabular}
\end{center}
\nota{Numbers of isolated exceptional fillings on $M_i$ up to the action of the symmetry group of $M_i$.}
\label{exceptional_2:table}
\end{table}

The numbers of isolated exceptional fillings on $M_i$ for $i=1,\ldots, 7$ considered up to the action of $\Iso(M_i)$ are listed in Table \ref{exceptional_2:table}. These numbers are smaller than those of Table \ref{exceptional_1:table}, but yet too big for our purposes, especially with the least symmetric and largest manifold $M_7$. We now want to reduce them further.

Recall that each $M_{i-1}$ is a filling of $M_i$. We use the standard meridian/longitude basis to identify the slopes on the boundary tori of $M_i$ with $\matQ \cup \{\infty\}$. If $s$ is a set of slopes we denote by $M_i(s)$ the manifold obtained by filling $M_i$ via $s$. Via a blow-down as in Figure \ref{Rolfsen:fig} we see that
$$M_{i-1} = M_i(-1)$$ 
for all $i\geq 2$. Note that there is no need of specifying which boundary component is filled thanks to the cyclic symmetry of all chain links. We now say that an exceptional filling of $M_i$ \emph{factors through $M_{i-1}$} if it contains $-1$ or any slope in the orbit of $-1$ along the action of $\Iso(M_i)$. Since we are classifying the exceptional slopes of $M_i$ inductively on $i$, it is natural to exclude those that factor through $M_{i-1}$. The surviving slopes are then collected in Table~\ref{exceptional_3:table}.

\begin{table}
\begin{center}
\begin{tabular}{lrrrrrrrr}
\hline
& \multicolumn{5}{c}{Number of cusps filled} & \\ \cline{2-8}
Manifold & 1 & 2 & 3 & 4 & 5 & 6 & 7 & Total\\ \hline
$M_1$ & 6    &     &     &    & & & & 6\\
$M_2$ & 6    & 4  &     &    & & & & 10\\
$M_3$ & 5    & 2  & 8  &  & & & & 15 \\
$M_4$ & 2    & 0  & 0   & 1  & & & & 3 \\
$M_5$ & 1    & 0  & 0  & 0  & 2  & & & 3 \\
$M_6$ &  1  &  0 &  0 &  2  &  4 &  40 & & 47 \\
$M_7$ & 2 & 1 & 2 & 6 & 61 & 313 & 1622 & 2007  \\
\hline
\alzapoco
\end{tabular}
\end{center}
\nota{Numbers of isolated exceptional fillings on $M_i$ that do not factor through $M_{i-1}$, up to the action of the symmetry group of $M_i$.}
\label{exceptional_3:table}
\end{table}

The numbers in Table \ref{exceptional_3:table} are extraordinarily small for $M_1,\ldots, M_5$ and are quite reasonable also for $M_6$. We can say informally that every $M_i$ adds a very small number of exceptional fillings to those of $M_{i-1}$ when $i\leq 5$. Only 6+10+15+3+3+47=84 basic exceptional fillings generate all the exceptional fillings of the manifolds $M_1,\ldots, M_6$. These 84 exceptional fillings are described in the tables at the end of the paper. The following theorem summarizes these discoveries.

\begin{teo}
The numbers of isolated exceptional fillings on $M_i$ that do not factor through $M_{i-1}$, up to the action of the symmetry group of $M_i$, are listed in Table \ref{exceptional_3:table}. The fillings of $M_1,\ldots, M_6$ are described in the Tables  \ref{M1_exc:table}, \ref{M2_exc:table}, \ref{M3_exc:table}, \ref{M4_exc:table}, \ref{M5_exc:table}, \ref{M6_exc:table}, \ref{M6b_exc:table}, \ref{M6c_exc:table}, and \ref{M6d_exc:table}.
\end{teo}

In the tables at the end of the paper we can find the exceptional slopes and a description of the 84 filled exceptional manifolds. Among these, 81 are graph manifolds and 3 are irreducible manifolds whose JSJ contains a hyperbolic piece, the figure-eight complement $M_1$. The manifold $M_6$ is the first manifold in the list that has some exceptional fillings that are not graph manifolds. A precise description of all the exceptional fillings of $M_6$ is given in Theorem \ref{6:teo}. 

\subsection{The manifold $M_7$} 
We are left with the 2,007 exceptional fillings of $M_7$ from Table \ref{exceptional_3:table}. These are yet too many to be reproduced here. Why does $M_7$ have such a big number of isolated fillings that do not factor through $M_6$? This is probably due again to its lack of symmetries: the group $\Iso(M_i)$ acts transitively on the cusps of $M_i$ for all $i$, and the groups $\Iso(M_4)$, $\Iso(M_5)$, and $\Iso(M_6)$ have a formidable amount of additional symmetries that send the slope $-1$ on any boundary torus $T$ to the sets of slopes (respectively) 
$$\left\{-1, \frac 12, \frac 32, 3\right\}, \qquad \left\{-1, \frac 12, 2\right\}, \qquad \{-1,1\}$$
on any boundary torus. Therefore $M_{i-1}$ is a filling of $M_i$ in multiple ways (respectively: in 16, 15, and 12 different ways), and hence there are many possibilities for an exceptional set of slopes on $M_i$ to factor through $M_{i-1}$. (We remark that the slopes $\{0,1,2,\infty\}$, $\{0,1,\infty\}$, and $\{0,\infty\}$ are exceptional on $M_4$, $M_5$, and $M_6$ respectively.)

On the other hand $\Iso(M_7)$ acts trivially on the slopes of a single boundary torus, so $M_6$ is a filling of $M_7$ in only 7 distinct ways. Since $\{0, \infty\}$ are exceptional, the first important non-exceptional slopes are $-1$ and $1$. It is natural to expect that most exceptional fillings of $M_7$ should contain either $-1$ or $1$. This is indeed the case, as Table \ref{exceptional_5:table} shows quite impressively. We denote by $N_6$ the hyperbolic manifold $N_6=M_7(1)$. 

\begin{table}
\begin{center}
\begin{tabular}{lrrrrrrrr}
\hline
& \multicolumn{5}{c}{Number of cusps filled} & \\ \cline{2-8}
Manifold & 1 & 2 & 3 & 4 & 5 & 6 & 7 & Total\\ \hline
$M_7$ & 2 & 0 & 0 & 0 & 1 & 14 & 73 & 90 \\
\hline
\alzapoco
\end{tabular}
\end{center}
\nota{Numbers of isolated exceptional fillings on $M_7$ that do not contain the slopes $-1$ and $1$, up to the action of the symmetry group of $M_7$.}
\label{exceptional_5:table}
\end{table}

\subsection{Another sequence of links}
\begin{figure}
\begin{center}
\includegraphics[width = 14 cm]{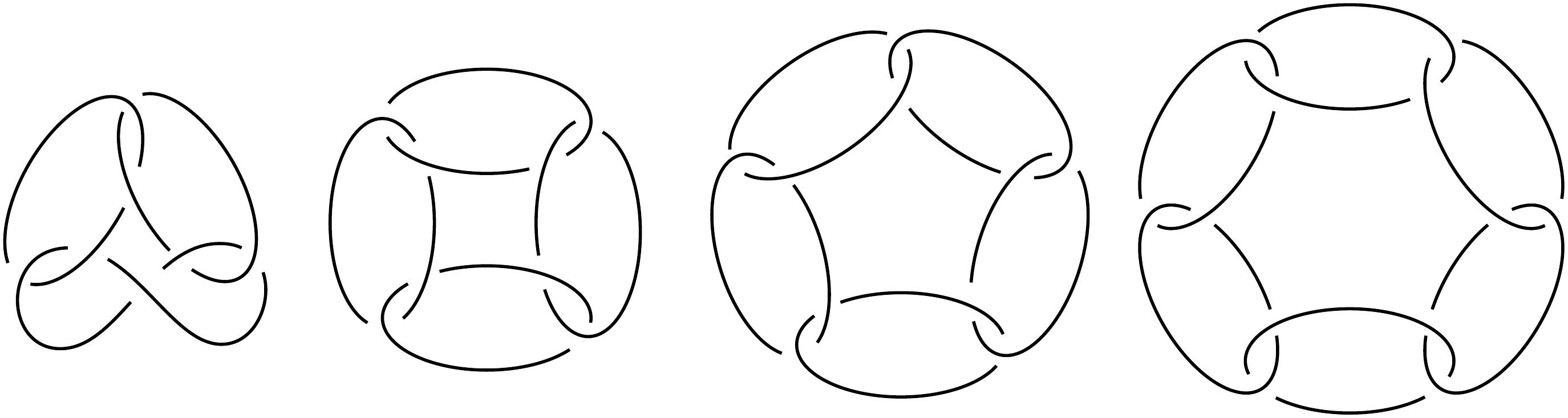}
\nota{Another sequence of hyperbolic chain links. Each link complement is a $1$-filling of the subsequent one.}
\label{sequence2:fig}
\end{center}
\end{figure}

We are still left with the problem of listing all the exceptional fillings of the new manifold $N_6 = M_7(1)$. By mirroring the blow-down in Figure \ref{Rolfsen:fig} we see that
$N_6$ is also the complement of a chain link with six components. It will be convenient to see $N_6$ as the last member of another sequence of chain link complements shown in Figure \ref{sequence2:fig}, that parallels somehow that of Figure \ref{new_sequence:fig}.

For every $i=3,\ldots, 6$, let $N_i$ be the complement of the chain link in Figure \ref{sequence2:fig} with $i$ components. These are all hyperbolic. Each $N_i$ is a $1$-filling of $N_{i+1}$, and $N_6$ is a $1$-filling of $M_7$. The manifolds $N_3, \ldots, N_5$ cannot be obtained as a Dehn filling of $M_6$ because they have some interesting exceptional fillings that $M_6$ does not have, as we will see. The volumes of the manifolds $M_i$ and $N_i$ are shown in Table \ref{volumes:table}.

\begin{table}
\begin{center}
\begin{tabular}{cccccc}
\hline
\alza $i$ & 3 & 4 & 5 & 6 & 7 \\
\hline
\alza ${\rm Vol}(M_i)$ & 5.33348 & 7.32772 & 10.14941 & 14.65544 & 19.79685 \\
\alza ${\rm Vol}(N_i)$ & 7.70691 & 10.14941 & 12.84485 & 16.00046 & \\
\hline
\end{tabular}
\end{center}
\alzapoco
\nota{The approximated volumes of the manifolds $M_i$ and $N_i$ for $i\geq 3$.}
\label{volumes:table}
\end{table}

It might be interesting to compare the numbers of exceptional fillings of the sequence $N_i$ with those of $M_i$. These are listed in Tables \ref{exceptional_N1:table} and \ref{exceptional_N2:table}. The symmetries of $N_3, \ldots, N_6$ are only those of the corresponding chain links, so they form a group of order 12, 16, 20, and 24.

\begin{table}
\begin{center}
\begin{tabular}{lrrrrrrrr}
\hline
& \multicolumn{5}{c}{Number of cusps filled} & \\ \cline{2-8}
Manifold & 1 & 2 & 3 & 4 & 5 & 6 & Total\\ \hline
$N_3$ & 6 & 48 & 70 & &  & & 124 \\
$N_4$ & 8 & 56 & 108 & 315 & & & 487 \\
$N_5$ & 10 & 50 & 155 & 695 & 3137 & & 4047  \\
$N_6$ & 12 & 33 & 150 & 1092 & 4962 & 28979 & 35228  \\
\hline
\alzapoco
\end{tabular}
\end{center}
\nota{Numbers of isolated exceptional fillings on $N_i$.}
\label{exceptional_N1:table}
\end{table}

\begin{table}
\begin{center}
\begin{tabular}{lrrrrrrrr}
\hline
& \multicolumn{5}{c}{Number of cusps filled} & \\ \cline{2-8}
Manifold & 1 & 2 & 3 & 4 & 5 & 6 & Total\\ \hline
$N_3$ & 2 & 10 & 14 & &  & & 26 \\
$N_4$ & 2 & 11 & 15 & 54 & & & 82 \\
$N_5$ & 2 & 7 & 19 & 71 & 326 & & 425 \\
$N_6$ & 2 & 5 & 14 & 98 & 418 & 2478 & 3015  \\
\hline
\alzapoco
\end{tabular}
\end{center}
\nota{Numbers of isolated exceptional fillings on $N_i$, up to the action of the symmetry group of $N_i$.}
\label{exceptional_N2:table}
\end{table}

\subsection{Some notable fillings}
Recall that our goal is to describe the exceptional fillings of $N_6$ with the minimum amount of information. 
Using SnapPy we discover some notable fillings of $N_3,\ldots, N_6$ in Table  \ref{some_fillings:table}. The table shows that the fillings $\{-3, -2, -1, 1\}$ on each $N_3, N_4, N_5$ and the fillings $\{-2, -1, 1\}$ on $N_6$ are diffeomorphic to either $N_{i-1}$ or some fillings of $M_3, M_4, M_5, M_6$. Since we have already examined the exceptional fillings of these manifolds, we disregard them: we say that a filling of $N_3, \ldots, N_6$ \emph{factors} if it contains one of these slopes (that is $\{-3,-2,-1,1\}$ for $N_3, N_4, N_5$, and $\{-2, -1, 1\}$ for $N_6$). We are happy with this definition because the number of isolated exceptional fillings of $N_3,\dots, N_6$ that do not factor is very small: see Table \ref{exceptional_4:table}. 

The following theorem summarizes our discoveries.

\begin{teo}
The numbers of isolated exceptional fillings of $N_i$ that do not factor, up to the action of the symmetry group of $N_i$, are listed in Table \ref{exceptional_4:table}. The fillings are described in the Tables \ref{N3:table}, \ref{N4:table}, \ref{N5:table}, and \ref{N6:table}.
\end{teo}

Among these 3+3+5+10 = 21 exceptional fillings, we find 14 graph manifolds and 7 irreducible manifolds whose JSJ contains a hyperbolic piece. The hyperbolic pieces that arise are $M_1, M_2, M_3$, and $M_4$.

\begin{table}
\begin{center}
\begin{tabular}{ccccc}
\hline
& \multicolumn{4}{c}{Slope} \\ \cline{2-5}
Manifold & $-3$ & $-2$ & $-1$ & 1 \\ \hline
$N_3$ &  $M_4(4, \frac 23)$   &  $M_3(-2)$  &  $M_3(-1) = M_2$   &  $M_6(-2, -\frac 12, \cdot, \frac 12, \frac 32)$ \\
$N_4$ &  $M_5(3, \frac 23)$   & $M_4(-2)$  & $M_4(-1) = M_3$    &  $N_3$ \\
$N_5$ & $M_6(2, 2)$    &  $M_5(-2)$ & $M_5(-1) = M_4$  & $N_4$  \\
$N_6$ &   & $M_6(-2)$ &  $M_6(-1) = M_5$   &  $N_5$   \\
\hline
\alzapoco
\end{tabular}
\end{center}
\nota{Some fillings of $N_i$ are diffeomorphic either to $N_{i-1}$ or to some filling of $M_3, M_4, M_5, M_6$. For instance from this table we infer that $N_4(-3) = M_5(3,\frac 23)$ and $N_5(-1) = M_5(-1) = M_4$.}
\label{some_fillings:table}
\end{table}

\begin{table}
\begin{center}
\begin{tabular}{lrrrrrrrr}
\hline
& \multicolumn{5}{c}{Number of cusps filled} & \\ \cline{2-8}
Manifold & 1 & 2 & 3 & 4 & 5 & 6 & 7 & Total\\ \hline
$N_3$ &  2   &  0   &  1   &    & & & & 3 \\
$N_4$ &  2   & 0  &  0   &  1  & & & & 3 \\
$N_5$ &  2   & 0  & 0  & 0  & 3 & & & 5 \\
$N_6$ &  2  & 0  &  0  & 0  & 2 & 6 & & 10 \\
\hline
\alzapoco
\end{tabular}
\end{center}
\nota{Numbers of isolated exceptional fillings on $N_i$ that do not factor, up to the action of the symmetry group of $N_i$.}
\label{exceptional_4:table}
\end{table}

\subsection{A final improvement}
We conclude this discussion by further reducing the numbers of Table \ref{exceptional_5:table}. Using SnapPy, we note the isometry $M_7(-2,-2) = N_6(-3)$. Since we have already classified the exceptional fillings of $N_6$, we may disregard all the fillings containing $(-2,-2)$. We say that a filling of $M_7$ \emph{factors} if it contains the slope $1$, $-1$, or the pair $(-2,-2)$ in two consecutive boundary components. The final survivors are counted in Table \ref{exceptional_6:table}. We will identify them in Tables \ref{M7a:table} and \ref{M7b:table}.

\begin{table}
\begin{center}
\begin{tabular}{lrrrrrrrr}
\hline
& \multicolumn{5}{c}{Number of cusps filled} & \\ \cline{2-8}
Manifold & 1 & 2 & 3 & 4 & 5 & 6 & 7 & Total\\ \hline
$M_7$ & 2 & 0 & 0 & 0 & 0 & 2 & 11 & 15 \\
\hline
\alzapoco
\end{tabular}
\end{center}
\nota{Numbers of isolated exceptional fillings on $M_7$ that do not contain the slopes $-1,1$, and $(-2,-2)$, up to the action of the symmetry group of $M_7$.}
\label{exceptional_6:table}
\end{table}

We summarize our discoveries on $M_7$.

\begin{teo}
The numbers of isolated exceptional fillings of $M_7$ that do not factor, up to the action of the symmetry group of $M_7$, are listed in Table \ref{exceptional_6:table}. The fillings are described in the Tables \ref{M7a:table} and \ref{M7b:table}.
\end{teo}

Among the exceptional fillings of $M_7$ we find infinitely many pairwise non-diffeomorphic closed manifolds whose JSJ has a hyperbolic piece. A complete description of all the manifolds that can be obtained as exceptional fillings of $M_7$ is given in Theorem \ref{except:M7:all:teo}.

\begin{rem} \label{correction:rem}
The tables in \cite{MaPeRo} of all the closed isolated exceptional fillings of $M_5$ contain
a few mistakes that we correct here. In \cite[Tables 9 and 10]{MaPeRo} the fillings $(-1, -2, -1, -3, -2)$ and $(-1, -2, -1, \frac 12, \frac 12)$ should be replaced with the correct ones $(-3, -2, -1, -3, -2)$ and $(-1, 2, -1, \frac 12, \frac 12)$ respectively.
On the other hand, the exceptional fillings $(-1, -\frac 12, -1, \frac 12, \frac 53)$, $(-1, -\frac 13, -1, \frac 23, \frac 32)$, $(-1, \frac 12, 3, -1, -\frac 12)$, and $(-1, \frac 12, -1, \frac 13, \frac 32)$ should be removed from \cite[Tables 9 and 10]{MaPeRo} because they are actually not isolated. For this reason the wrong numbers $5232$ and $52$ appeared in \cite[Tables 1 and 2]{MaPeRo} instead of the correct ones $4818$ and $48$ that we display here. 
\end{rem}

\section{The exceptional fillings}

In the previous section we have reduced all the isolated exceptional fillings of the hyperbolic manifolds $M_1,\ldots, M_7$ to some, as small as possible, ``generating'' set. We now describe explicitly these generating exceptional fillings. 

\subsection{Notation}
We use the same notation of \cite{Ma, MaPeRo} for Seifert and graph manifolds, which seems standard. We quickly recall it here. Given a compact surface $\Sigma$, possibly with boundary, and some pairs $(p_1,q_1), \ldots, (p_k,q_k)$ of coprime integers, the notation
$$X = \big( \Sigma, (p_1,q_1), \ldots, (p_k,q_k) \big)$$
denotes the 3-manifold $X$ obtained as follows. We remove $k$ open discs from $\Sigma$, thus getting a new surface $\Sigma'$. Then we attach $k$ solid tori to the (unique) oriented circle bundle over $\Sigma'$ by killing the slopes $(p_1,q_1), \ldots, (p_k,q_k)$ in any $k$ boundary tori. We use here as a basis a meridian in $\partial \Sigma'$ and a longitude $\{{\rm pt}\} \times S^1$, oriented as a positive basis.

Note that the case $p_i = 0$ is allowed here. It is a standard fact on Seifert manifolds that if $p_i\neq 0$ for all $i$ then $X$ is a Seifert manifold, while if $p_i=0$ for some $i$ then $X$ ``degenerates'' to a connected sum of lens spaces and solid tori. More specifically, if $\Sigma$ is orientable we have
$$\big( \Sigma, (0,1), (p_2,q_2), \ldots, (p_k,q_k) \big)
=
\lens {p_2}{q_2} \# \cdots \# \lens{p_k}{q_k} \#_{2g} (S^2 \times S^1) \#_b(D^2 \times S^1)
$$
where $g$ and $b$ are the genus and the number of boundary components of $\Sigma$.
Concretely, in most cases the surface $\Sigma$ will be either $S^2, D, A$, or $P$, that is a sphere, a disc, an annulus, or a pair of pants. 

The manifold $X$ described in this way is naturally equipped with an orientation and a standard meridian/longitude basis on each boundary torus. There is no need to distinguish between the boundary tori of a Seifert manifold since they are all equivalent up to diffeomorphism. With that in mind, given some manifolds $X, Y, Z,$ and matrices $A, B \in \GL(2,\matZ)$ we may write
$$X \bigcup\nolimits_{A} Y \bigcup\nolimits_B Z, \qquad X/_A, \qquad X \bigcup\nolimits_A\nolimits^B Y$$
to denote some graph manifolds that decomposes along tori into the pieces $X,Y,Z$, glued via the maps $A, B$. In the second example two distinct boundary components of $X$ are identified via $A$. In the third, two manifolds $X, Y$ have two pairs of boundary tori glued via $A$ and $B$. All the matrices here will have $\det = -1$. We also use the notation $T_A$
to denote a torus fibration over $S^1$ with monodromy $A$, and in this case we have $\det A = 1$.

\subsection{Ambiguities} \label{ambiguities:subsection}
The same graph manifold may be described in various different ways and unfortunately in many occasions there is no preferred description. 

A useful set of moves that modify the notation of a graph manifold was collected in \cite[Lemma 2.1]{MaPe}. We report here for completeness the ones that are more relevant for us: we will use them at various points. Here are the first ones:
\begin{eqnarray}
\big( \Sigma, (a,b), \ldots (y,z) \big) & = & \big( \Sigma, (a,-b), \ldots (y,-z) \big)
    \label{minus:eqn} \\
\seifdueeul \Sigma abcd \ldots & = & \seifdueeul \Sigma{a}{b+ka}{c}{d-kc}\ldots \label{from:D2:to:D2:eqn}\\
\seifunoeul \Sigma ab \ldots & = & \seifunoeul \Sigma a{b+ka} \ldots 
    \quad {\rm if\ }\partial \Sigma \neq \emptyset \label{from:F:to:F:eqn}\\
\seifdueeul \Sigma10ab\ldots & = & \seifunoeul {\Sigma}ab\ldots \label{from:D2:to:D1:eqn}
\end{eqnarray}

In the first move we change the orientation of the manifold.  
In this paper we have chosen to write a Seifert manifold as standardly as possible, with positive normalized numbers if the manifold has boundary: the unique oriented Seifert manifold fibering over the orbifold $(D,2,2)$ is usually denoted as $\seifdue D 2121$, although the notation $\seifdue D212{-1}$ would also be natural since it visibly expresses the fact that the Euler number vanishes; there are two oriented Seifert manifolds fibering over $(D,2,3)$, and these are $\seifdue D2131$ and $\seifdue D2132$. They are orientation-reversingly diffeomorphic.

The following moves involve the gluing of two pieces:
\begin{eqnarray}
X \bigcup\nolimits_A Y & = & X \bigcup\nolimits_{-A} Y \label{meno:eqn} \\
\seifunoeul \Sigma{a}{b}\ldots\bigu mnpq X & = & \seifunoeul \Sigma{a}{b+k a}\ldots\bigu
    {m+k n}n{p+k q}q X \label{from:D2X:to:D2X:eqn}\\
X\bigu mnpq \seifunoeul \Sigma ab\ldots & = & X\bigu mn{p-k m}{q-k n}
    \seifunoeul \Sigma a{b+ka}{\ldots} \label{from:XD2:to:XD2:eqn} 
\end{eqnarray}

The moves (\ref{from:D2X:to:D2X:eqn}, \ref{from:XD2:to:XD2:eqn}) also apply when two boundary tori of the same block are glued together, but move (\ref{meno:eqn}) does not! Note that with (\ref{meno:eqn}, \ref{from:D2X:to:D2X:eqn}, \ref{from:XD2:to:XD2:eqn}) it is not possible to change the absolute value $|n|$ of the top-right element $n$ in the gluing matrix. Indeed $|n|$ has an important geometric significance: it is the geometric intersection of the fibers of the two glued Seifert manifolds.

There are also some more complicated moves that occur in more sporadic cases. The following reflect the fact that $\seifdue D2121$ has an alternative description as the orientable circle bundle $S\timtil S^1$ over the M\"obius strip $S$.
\begin{eqnarray}
\seifdue D2121 \bigu mnpq  X & = &
(S\timtil S^1)\bigu n{n-m}q{q-p} X
    \label{from:D2X:to:SX:eqn} \\
X\bigu mnpq \seifdue D2121 & = & X \bigu{m+p}{n+q}{-m}{-n} (S\timtil S^1)
    \label{from:XD2:to:XS:eqn} 
\end{eqnarray}

Finally, it is sometimes useful to understand how the manifold ``degenerates'' when we perform a longitudinal filling:
\begin{eqnarray}
\seiftre {S^2}abcd 01 & = & \lens ab\# \lens cd
    \label{from:D2:to:lens:plus:lens:eqn}  \\
\seifdue D01ab\bigu mnpq \seifeul \Sigma\ldots & = & \lens ab\#\seifunoeul {\Sigma'}nq\ldots
    \label{from:DX:to:lensX:eqn}
\end{eqnarray}

Here $\Sigma'$ is $\Sigma$ with one boundary component capped off. In this paper $S^2 \times S^1$ is also denoted as the lens space $L(0,1)$.

\subsection{Zero and infinity}
Let $M$ be the complement of any chain link $L\subset S^3$. The fillings $0$ and $\infty$ on any boundary component of $M$ are easily understood, and they are always exceptional.

\begin{figure}
\begin{center}
\includegraphics[width = 16 cm]{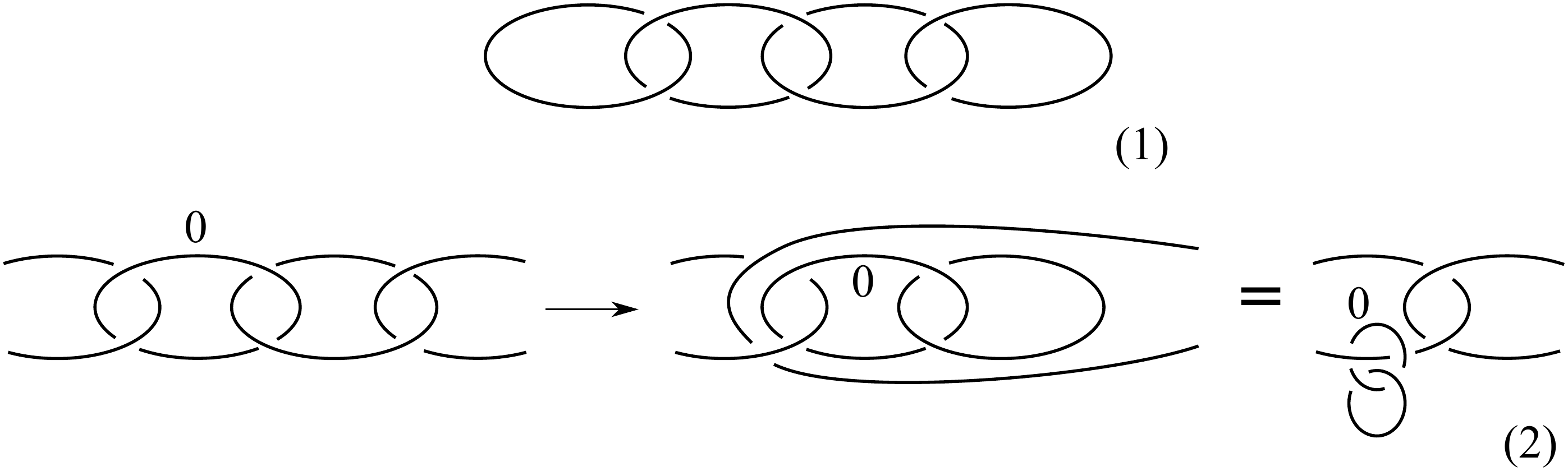}
\nota{The fillings $\infty$ and $0$ on any chain link.}
\label{open_chain:fig}
\end{center}
\end{figure}

The filling $\infty$ corresponds to the removal of a component from $L$, so we get the complement of an open chain in $S^3$ as in Figure \ref{open_chain:fig}-(1), that is easily identified as the graph manifold $A \times S^1, P \times S^1$, or
$$P\times S^1 \bigu 0110 \cdots \bigu 0110 P \times S^1$$
depending on the number of components of $L$.

The filling $0$ may be modified with a handle slide as shown in Figure \ref{open_chain:fig}-(2). The resulting manifold is the complement of another chain link (with two components less) attached to a $P \times S^1$.

These exceptional fillings $0$ and $\infty$ will appear on all the tables concerning  the various chain links studied in this paper.

\subsection{The exceptional fillings}
We can finally describe the exceptional fillings of $M_6$ and $M_7$. Let us start with $M_6$. For completeness, we also review all the exceptional fillings of $M_1,\ldots, M_5$, already classified in \cite{MaPe, MaPeRo}. All the tables are postponed to the end of the paper for the sake of clarity.

\begin{teo} \label{exc:6:teo}
Every isolated exceptional filling of $M_i$ with $i=1,\ldots, 6$ that does not factor through $M_{i-1}$ is equivalent, up to the action of $\Iso(M_i)$, to precisely one of those listed in Tables \ref{M1_exc:table}, \ref{M2_exc:table}, \ref{M3_exc:table}, \ref{M4_exc:table}, \ref{M5_exc:table}, \ref{M6_exc:table}, \ref{M6b_exc:table} and \ref{M6c_exc:table}. The filled manifold is also shown there.
\end{teo}

Recall that factoring through $M_{i-1}$ is equivalent to containing $-1$ in some cusp, or any other slope obtained from $-1$ by the action of $\Iso(M_i)$. The tables show the isolated exceptional slopes (one representative for each orbit of $\Iso(M_i)$), the filled manifold, and its integral first homology group. The notation for the filled manifold sometimes differ from \cite{MaPe, MaPeRo} via some of the moves listed in Section \ref{ambiguities:subsection}.

We now do the same with $N_3,\ldots, N_6$. Recall that ``factoring'' here means that the filling slopes contain at least one of the numbers $\{-3,-2,-1,1\}$ for $N_3, N_4, N_5$ and of $\{-2,-1,1\}$ for $N_6$.

\begin{teo} \label{exc:N:teo}
Every isolated exceptional filling of $N_3, \ldots, N_6$ that does not factor is equivalent, up to the action of $\Iso(N_i)$, to precisely one of those listed in Tables 
\ref{N3:table}, \ref{N4:table}, \ref{N5:table}, and \ref{N6:table}.
The filled manifold is also listed there.
\end{teo}

Finally, we turn to $M_7$. Recall that ``factoring'' here means that the filling slope contains either $-1$, $1$, or the pair $(-2,-2)$ in two consecutive boundary tori.

\begin{teo} \label{exc:7:teo}
Every isolated exceptional filling of $M_7$ that does not factor is equivalent, up to the action of $\Iso(M_7)$, to precisely one of those listed in Tables \ref{M7a:table} and \ref{M7b:table}.
\end{teo}

The tables shown so far contain a fair amount of information. From these, we can easily deduce which kinds of non-hyperbolic filling we can obtain from each manifold $M_i$. We do this in the following sections.

\subsection{The manifold $M_5$}

The following theorem was already proved in \cite[Corollary 1.3]{MaPeRo}. 

\begin{teo} \label{5:teo}
The closed non-hyperbolic fillings of $M_5$ are precisely the manifolds:
$$
\seifdue D abcd \bigu 0110 \seifdue D efgh,
$$
$$\seifdue{D}2121 \bigu {1+n}{2+n}{-n}{-1-n} \seifdue{D}2131,$$
$$\seifuno{A}ab \bigb 0110, \qquad \seifuno A 21 \bigb 120{-1}$$
where $(a,b)$, $(c,d)$, $(e,f)$, $(g,h)$ are arbitrary pairs of coprime integers, and $n\in \{0,1,2,3\}$.
\end{teo}

We note that the first family in the theorem contains many different kinds of manifolds.
\begin{prop} \label{a:prop}
The manifolds $X$ that may be obtained via the description
$$X = \seifdue D abcd \bigu 0110 \seifdue D efgh$$
are precisely the following:
\begin{enumerate}
\item The lens spaces and connected sums of two lens spaces.
\item The Seifert manifolds fibering over $S^2$ with 3 exceptional fibres.
\item The Seifert manifolds fibering over $\matP^2$ with 2 exceptional fibres.
\item The Seifert manifold $(K,1)$.
\item The graph manifolds whose JSJ decomposition is as in the description of $X$.
\end{enumerate}
\end{prop}

Here $K$ is the Klein bottle and $(K,1)$ is the fibration over $K$ with Euler number 1.

\begin{proof}
We use the moves described in Section \ref{ambiguities:subsection}. If one of $a,c,e,g$ is zero, we get a connected sum of two lens spaces. Now suppose $a,c,e,g \neq 0$. If $a=1$, we get a Dehn filling of $\seifdue Defgh$, hence either a lens space, a connected sum of two lens spaces, or a Seifert manifold fibering over $S^2$ with 3 exceptional fibres.

We are left with the case $|a|, |c|, |e|, |g| \geq 2$. In general, we get a graph manifold whose JSJ decomposition is as in the description of $X$. There is only one exceptional case to consider: if $a=c=2$, then up to some moves we get
\begin{align*}
X & = \seifdue D 2121 \bigu k110 \seifdue D efgh \\
&= \big( S \timtil S^1 \big) \bigu 1{1-k}0{-1}  \seifdue D efgh
\end{align*}
for some $k\in \matZ$. If $k=1$ the fibers of the two blocks match to give a Seifert manifold with two exceptional fibres over $\matRP^2$. If $e=g=2$ also the right block has another fibration and we get
$$X = \big(S \timtil S^1) \bigu 1{1-k} l{l-1-kl} \seifdue D2121 = 
\big(S \timtil S^1 \big) \bigu {1+l}{l-k-kl}{-1}{k-1}  \big(S \timtil S^1 \big)$$
for some $k,l\in \matZ$.
The two fibres match when $l-k-kl=0$ and $k-1=\pm1$. We get two cases $(k,l)=(0,0)$ or $(2,-2)$ and in both cases we get the Seifert manifold $(K,\pm 1)$.
\end{proof}

The generic case (5) consists precisely of all the irreducible 3-manifolds whose JSJ decomposition consists of two Seifert pieces, each fibering over a disc with two cone points, whose fibers meet in the glued torus with geometric intersection one. Theorem \ref{5:teo} says that $M_5$ has also some more sporadic exceptional fillings where this geometric intersection is 2, 3, 4, or 5.

Among the exceptional fillings of $M_5$ we also find another family that can be analyzed in a similar fashion:
\begin{prop} \label{b:prop}
The manifolds $X$ that may be obtained via the description
$$\seifuno{A}ab \bigb 0110$$
are precisely the following:
\begin{enumerate}
\item The manifold $S^2\times S^1$.
\item The torus bundles of type
$$T_{\tiny{\matr b{1}{-1}0}}.$$
\item The graph manifolds whose JSJ decomposition is as in the description of $X$.
\end{enumerate}
\end{prop}
\begin{proof}
When $a=0$ we get $S^2\times S^1$. When $a=1$ we get the torus bundle
$$T_{\tiny{\matr b{-1}{1}0}} = T_{\tiny{\matr b{1}{-1}0}}.$$
When $|a| \geq 2$ we get a manifold whose JSJ decomposition is as described. 
\end{proof}

As above, the manifolds that we get in (3) are precisely all the irreducible 3-manifolds whose JSJ decomposition consists of a single piece fibering over an annulus with a cone point, whose fibers meet in the glued torus with geometric intersection one. Theorem \ref{5:teo} exhibits also a sporadic example with geometric intersection 2.

As we already knew from \cite{MaPeRo}, all the exceptional fillings of $M_5$ are graph 
manifolds. We now discover here that this is not the case for $M_6$.

\subsection{The manifold $M_6$.}
We now turn to the exceptional fillings of $M_6$.

\begin{teo} \label{6:teo}
The closed non-hyperbolic fillings of $M_6$ are precisely the manifolds:
$$\seifdue Dabcd \bigu 0110 \seifuno Aef \bigu 0110 \seifdue Dghij,$$
$$\seifuno Aab
\bigcup\nolimits_{{\tiny{\matr 0110}}}^{{\tiny{\matr 0110}}} \seifuno Acd,$$ 
$$M_1 \bigu {-1}011 \seifdue D2121, \qquad M_1 \bigu {-1}110 \seifdue D2121$$
where $(a,b)$, $(c,d)$, $(e,f)$, $(g,h)$, $(i,j)$ are arbitrary pairs of coprime integers.
\end{teo}
\begin{proof}
The Tables \ref{M6_exc:table}, \ref{M6b_exc:table}, \ref{M6c_exc:table} and \ref{M6d_exc:table} show that the exceptional fillings 
$$\infty, \quad \left(-2,-\dfrac{1}{2},\boldsymbol\cdot,\dfrac{1}{2},2\right), 
\quad \left(-3,-2,-\dfrac{1}{3},3,2,\dfrac{1}{3}\right), \quad
\left(-3,-\dfrac{3}{2},-\dfrac{1}{2},2,2,\dfrac{1}{3}\right)$$
give rise precisely to all the manifolds listed in the theorem. Conversely, a case by case analysis of the manifolds listed in Theorem \ref{5:teo} and Tables \ref{M6_exc:table}, \ref{M6b_exc:table}, \ref{M6c_exc:table}, and \ref{M6d_exc:table} shows that all the exceptional fillings of $M_6$ are of one of these types. Here are the details. Using the moves of Section \ref{ambiguities:subsection} we see that 
$$\seifdue Dabcd \bigu 0110 \seifuno A1f \bigu 0110 \seifdue Dghij
= $$
$$\seifdue Dabcd \bigu 1f0{-1} \seifdue Dghij,
$$
$$\seifuno Aab
\bigcup\nolimits_{{\tiny{\matr 0110}}}^{{\tiny{\matr 0110}}} \seifuno A1d
=  \seifuno Aab \bigb 1d0{-1}.
$$
By applying the moves (\ref{from:D2X:to:D2X:eqn},
\ref{from:XD2:to:XD2:eqn}) we deduce that we can actually obtain in this way all the manifolds of the following two types:
$$\seifdue Dabcd \bigcup\nolimits_B \seifdue Dghij, \qquad
\seifuno A ab /_B$$
where $B$ is any matrix that can be written as 
$$B = \begin{pmatrix} 
1 + mf & f \\
-m-n-mnf & -(1 + nf) 
\end{pmatrix}$$
for some $m,n,f\in \matZ$. 
In other words, here $B$ is any matrix $B = \tiny\matr rfst$ with $\det B = -1$ such that $r \equiv 1 \mod f$ and $t \equiv -1 \mod f$. When the manifold is of the first type we can also exchange $B$ with $-B$ using the move (\ref{meno:eqn}) from Section \ref{ambiguities:subsection}, so we may also get $r \equiv -1 \mod f$ and $t \equiv 1 \mod f$ in that case. 
All the manifolds that arise from Theorem \ref{5:teo} and Tables \ref{M6_exc:table}, \ref{M6b_exc:table}, \ref{M6c_exc:table} and \ref{M6d_exc:table} are of this kind, except of course the two
manifolds whose JSJ decomposition contains a hyperbolic piece.
\end{proof}

Theorem \ref{6:teo} exhibits a couple of important differences between $M_6$ and the manifolds $M_1, \ldots, M_5$. The first is that all the graph manifolds come into two big families, and there are no sporadic manifolds outside of these. The second is of course the presence of two irreducible manifolds whose JSJ decomposition contains some hyperbolic piece.

The following proposition furnishes some details on the graph manifolds produced by the first family.

\begin{prop}
The manifolds $X$ that may be obtained via the description
$$X = \seifdue D abcd \bigu 0110 \seifuno A ef \bigu 0110 \seifdue D ghij$$
are precisely the following:
\begin{enumerate}
\item The manifolds that arise in Proposition \ref{a:prop}.
\item The connected sums of three lens spaces.
\item The connected sums of a Seifert manifold over $S^2$ with 3 exceptional fibres and a lens space.
\item The Seifert manifolds over $S^2$ with 4 exceptional fibres.
\item The Seifert manifolds over $K$ with 0 or 1 exceptional fibres.
\item The graph manifolds whose JSJ decomposition is
$$X = \seifdue D abcd \bigcup\nolimits_B \seifdue D ghij$$
with
$$B = \begin{pmatrix} 
1 + mf & f \\
-m-n-mnf & -(1 + nf) 
\end{pmatrix}.$$
\item The graph manifolds whose JSJ decomposition is
$$X = \seifuno S ab \bigu 0110 \seifdue D cdef$$
\item The graph manifolds whose JSJ decomposition is as in the description of $X$.
\end{enumerate}
\end{prop}
Here $K$ and $S$ are the Klein bottle and the M\"obius strip.
\begin{proof}
If $a=0$ we get a connected sum $\lens cd \# \big(S^2, (f,e),(g,h),(i,j)\big)$. The second addendum may in turn give rise to a connected sum of two lens spaces. If $e=0$ we get a connected sum of two lens spaces. So we suppose that $a,c,e,g,i \neq 0$.

If $a=1$ we get a manifold as in Proposition \ref{a:prop}. So we suppose $|a|, |c|,  |g|, |i| \geq 2$. If $e=1$ we get a manifold
$$\seifdue Dabcd \bigcup\nolimits_B \seifdue Dghij$$
where $B$ is any matrix that can be written as 
$$B = \begin{pmatrix} 
1 + mf & f \\
-m-n-mnf & - (1 + nf) 
\end{pmatrix}.$$
See the proof of Theorem \ref{6:teo}. If $f=0$ we get (4). If $f\neq 0$ we get a graph manifold as in (6), except possibly when one (or both) piece is $\seifdue D 2121$ and the fibrations match: in this way we could obtain a Seifert fibration over $\matRP^2$ with two singular fibres or over $K$ without singular fibres; the former was already obtained in (1) and the latter will be obtained in the next paragraph by other means, so we ignore it. 

We can now suppose that also $|e| \geq 2$. We get a manifold of type (8), except when the left (or the right) block is $\seifdue D 2121$ and its alternative fibration matches with that of the central block. This may happen and we get a manifold of type (7), unless this happens on both extreme blocks simultaneously and in this case we get (5).
\end{proof}

We do the same analysis with the second family of graph manifolds.

\begin{prop}
The manifolds $X$ that may be obtained via the description
$$\seifuno Aab
\bigcup\nolimits_{{\tiny{\matr 0110}}}^{{\tiny{\matr 0110}}} \seifuno Acd$$ 
are precisely the following:
\begin{enumerate}
\item The manifolds $(S^2\times S^1) \# L(p,q)$.
\item The torus bundles $T_C$ with
$$C = \begin{pmatrix} 
-(1 + mf) & -f \\
-m-n-mnf &- (1 + nf) 
\end{pmatrix}.
$$
\item The graph manifolds whose JSJ decomposition is
$$X = \seifuno A ab /_B$$
with
$$B = \begin{pmatrix} 
1 + mf & f \\
-m-n-mnf &-(1 + nf) 
\end{pmatrix}.$$
\item All graph manifolds whose JSJ decomposition is as in the description of $X$.
\end{enumerate}
\end{prop}
\begin{proof}
If $a=0$ we get $(S^2 \times S^1) \# L(d,c)$. So we suppose $a, c \neq 0$. If $c=1$ we get (3), unless $a=1$ and in this case we get (2). If $|a|, |c| \geq 2$ we get (4).
\end{proof}

We note in particular that we get all the torus bundles with monodromy $\tiny \matr {-1}0a{-1}$. However, we do not get the identity matrix! We deduce the following.

\begin{cor}
Among the six orientable flat 3-manifolds, five can be obtained by Dehn filling $M_6$, but the 3-torus cannot.
\end{cor}
\begin{proof}
The four that fiber over $S^2$ with three exceptional fibers or over $\matRP^2$ with two exceptional fibers can already be obtained from the magic manifold $M_3$. The one that fibers over $S^2$ with four exceptional fibers (equivalently, over $K$) is obtained with $M_6$.
\end{proof}

One important novelty in the exceptional fillings of $M_6$ is of course the presence of two sporadic irreducible manifolds whose JSJ decomposition contains a hyperbolic piece. Both exceptional manifolds decompose into the figure-eight knot complement $M_1$ and the Seifert manifold $\seifdue D 2121$, which is diffeomorphic to the $I$-bundle $K\timtil I$ over the Klein bottle $K$ and to the orientable $S^1$-bundle $S \timtil S^1$ over the M\"obius stip $S$. By looking at the gluing matrices, we note that in both cases the meridian of the figure-eight complement (which is also the shortest curve in a flat cusp section) is attached to the fiber of the alternative fibration $S \timtil S^1$, that is represented as the slope $(-1, 1)$ in the fibration $\big(D, (2,1), (2,1) \big)$.

\subsection{The manifold $M_7$.} \label{M7:subsection}
We now turn to $M_7$. We recall that the simplest method we found to describe all the exceptional fillings of $M_7$ was to study an alternative sequence of chain links $N_3, \ldots, N_6$.

We first note the quite surprising fact that $N_3$ contains infinitely many distinct exceptional fillings with hyperbolic pieces. As an aside, this implies that the manifolds $N_3, N_4, N_5$ are not fillings of $M_6$. 

We can now list all the exceptional fillings of $M_7$.

\begin{teo} \label{except:M7:all:teo}
The closed non-hyperbolic fillings of $M_7$ are precisely the manifolds:
$$\seifdue Dabcd \bigu 0110 \seifuno Aef \bigu 0110 \seifuno Agh \bigu 0110 \seifdue Dijkl,$$
$$\seifuno Aab
\bigcup\nolimits_{{\tiny{\matr 0110}}}^{{\tiny{\matr 0110}}} \seifuno Acd,$$ 
$$ M_5\big((a,b), (c,d), (e,f), (g,h)\big) \bigu 0110 \seifdue Dijkl, $$
$$M_2(a,b) \bigu {-1}0{1}{1} \seifdue D 2121,$$
$$\seifuno A 2{1} \bigb {n-1}n11$$
where $(a,b)$, $(c,d)$, $(e,f)$, $(g,h)$, $(i,j)$, $(k,l)$ are arbitrary pairs of coprime integers and $n\in \{3, 4, 5, 6\}$. In the third family we suppose $|i|, |k| \geq 2$.
\end{teo}

We note the reappearance of some sporadic graph manifolds in the list, which were absent in $M_6$. The four sporadic graph manifolds listed in the last row are not members of the previous families.

\begin{proof}
We have to check that all the exceptional fillings of $M_7$ are of this type, and to this purpose we only need to verify this for the manifolds listed in Tables \ref{N3:table}, \ref{N4:table}, \ref{N5:table}, \ref{N6:table}, \ref{M7a:table}, and \ref{M7b:table}. Concerning graph manifolds, this is easily settled using the moves described in Section \ref{ambiguities:subsection} when necessary.

We are left with the non-graph manifolds with non-trivial JSJ decomposition. The manifolds $N_i(0)$ with $i=3,\ldots, 6$ are obviously a filling of $N_7(0)$ and hence can be excluded, since their fillings are already contained in the third family. The tables contain three manifolds 
$$X_i = M_2 \bigcup\nolimits_{A_i} \seifdue D 2121,$$
where $A_i$ is one of the matrices
$$A_1= \begin{pmatrix} -1 & 2 \\ 1 & -1 \end{pmatrix}, \qquad
A_2= \begin{pmatrix} -1 & 0 \\ 1 & 1 \end{pmatrix}, \qquad
A_3 = \begin{pmatrix} -1 & 1 \\ 1 & 0 \end{pmatrix}.
$$
There are also more manifolds where $M_2$ is replaced either by $M_1$ or by $M_2(-2)$ and $A_i$ is still of one of these three types. Since $M_1 = M_2(-1)$ via an isometry that acts as the identity on the other boundary torus, these manifolds are fillings of the $X_i$ and can be ignored. 
To conclude we need to show that the manifolds $X_1$ and $X_3$ are fillings of 
$$X = M_5 \bigu 0110 P \times S^1.$$
It will be convenient to use the moves in Section \ref{ambiguities:subsection} and write them as
\begin{align*}
X_1 & = M_2 \bigu {-1}201 \seifdue D 212{-1}, \\
X_3 & = M_2 \bigu {-1}1{-1}2 \seifdue D 2{-1}2{-1}.
\end{align*}

Using SnapPy we find an isometry from $M_2$ to $M_5(-1,-2,-2)$ that acts on a boundary torus as the matrix
$$B = \begin{pmatrix}
1 & -1 \\ 0 & 1 \end{pmatrix}.$$
We also note that there are isometries of $M_5$ that act on the cusps like the matrices
$$C_1 = \begin{pmatrix}
0 & 1 \\
-1 & 1
\end{pmatrix}, \qquad
C_2 = \begin{pmatrix}
-1 & 1 \\ 
-1 & 0
\end{pmatrix}
$$
We deduce that both $X_1$ and $X_3$ are Dehn fillings of $X$ because
$$\begin{pmatrix} -1 & 2 \\ 0 & 1 \end{pmatrix} = 
\begin{pmatrix} 0 & 1 \\ 1 & 0 \end{pmatrix} C_1B,
\qquad
\begin{pmatrix} -1 & 1 \\ -1 & 2 \end{pmatrix} =
\begin{pmatrix} 0 & 1 \\ 1 & 0 \end{pmatrix} C_2B.
$$
The proof is complete.
\end{proof}
\begin{cor}
The 3-torus $S^1\times S^1 \times S^1$ is not a filling of $M_7$.
\end{cor}
In particular, the Borromean rings complement is not a filling of $M_7$. As pointed out by the referee, this corollary is actually well-known: since the minimally twisted chain links are strongly invertible, every Dehn surgery on $M_7$ is a double branched cover over some link in $S^3$. By a result of Fox \cite{Fox}, the 3-torus cannot be realized as a double branched cover on any link in $S^3$. Therefore it cannot be realized as a Dehn surgery on any minimally twisted chain link.

\section{The cusped census}
Many cusped hyperbolic manifolds of the Callahan -- Hildebrand -- Weeks census \cite{CHW} are Dehn fillings of $M_7$ and a classification of their exceptional fillings can be deduced from the theorems stated here.

The Callahan -- Hildebrand -- Weeks census contains all the cusped hyperbolic manifolds $N$ with complexity $c\leq 7$ (the complexity is the minimum number of tetrahedra in a topological ideal triangulation). To see whether $N$ can be obtained as a filling of $M_i$ for some $i\leq 7$, a very simple and concrete method consists of drilling the shortest simple closed curve found by SnapPy in $N$ multiple times, until we get a manifold with $i$ cusps. This typically leads to some hyperbolic manifold $N_i$ and one checks via SnapPy whether $N_i$ is isometric to $M_i$ or not. Note that as soon as $N_i=M_i$ we also have $N_{i+1} = M_{i+1}$ because each $M_{i+1}$ is obtained by drilling a shortest curve from $M_i$. If this crude algorithm fails, of course this does not imply that $N$ is not a filling of $M_i$.

\begin{table}
\begin{center}
\begin{tabular}{lrrrrrrr}
\hline
$c$ & $M_3$ & $M_4$ & $M_5$ & $M_6$ & $M_7$ & left & Total
\\ \hline
$\leq 5$ & 216 & 42 & 8 & 25 & 0 & 10 & 301  \\
6 & 426 & 286 & 94 & 111 & 8 & 37 & 962 \\
7 & 1077 & 1142 & 558 & 519 & 105 & 151 & 3552 \\
\hline
Total $\leq 7$
& 1719 & 1470 & 660 & 655 & 113 & 198 & 4815 \\
\hline
\alzapoco
\end{tabular}
\end{center}
\nota{Numbers of cusped hyperbolic manifolds of complexity $c$ that transform into $M_i$ after repeatedly drilling the shortest curve found by SnapPy. If $i\geq 4$ we only count those that were not previously transformed into $M_{i-1}$.}
\label{census:table}
\end{table}

We display in Table \ref{census:table} the number of manifolds in each complexity $c\leq 5$, $c=6$, and $c= 7$ for which this algorithm produces a positive answer for $M_i$. We could represent all the manifolds of the census as a filling of $M_7$, except a number of 10, 37, and 151 of them with $c\leq 5, c=6$, and $c=7$ respectively.

\begin{figure}
\begin{center}
\includegraphics[width = 4 cm]{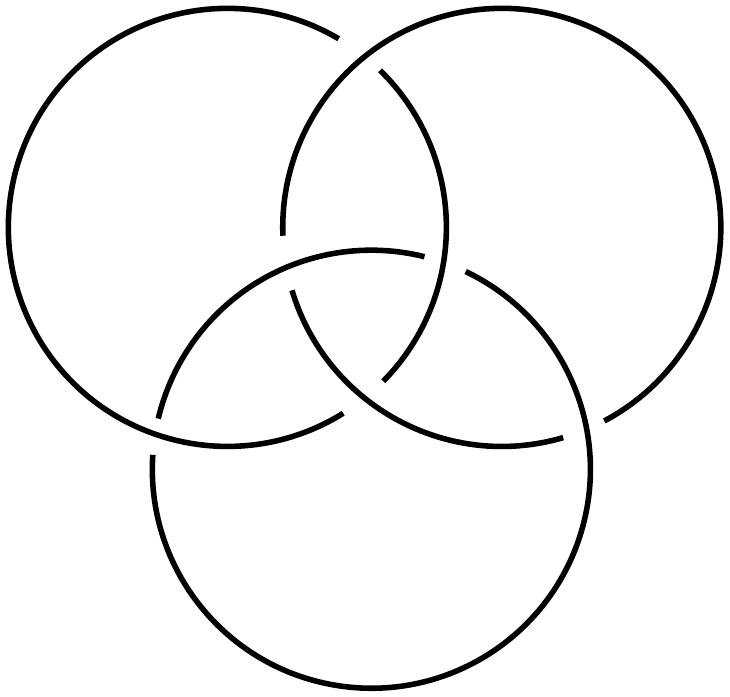}
\nota{The Borromean rings.}
\label{Borromean:fig}
\end{center}
\end{figure}

\subsection{The Borromean rings}
We have discovered at the end of Section \ref{M7:subsection} that the complement $W$ of the Borromean rings (shown in Figure \ref{Borromean:fig}) is not a Dehn filling of $M_7$. While performing the drilling algorithm, we notice that in each of the 10 remaining manifolds $N$ with $c(N)\leq 5$ of Table \ref{census:table}, the drilled manifold $N_3$ is in fact isometric to the Borromean rings complement $W$. More generally, among the 10, 47, and 198 remaining manifolds with $c \leq 5$, $c=6$, and $c=7$, for 10, 34, and 92 of them the manifold $N_3$ is isometric to $W$.

We are now led naturally to the study of the exceptional fillings of the Borromean ring complement $W$. These are easily classified. We note that $W(-1) = M_2$ is the Whitehead link complement, so we consider only the Dehn fillings that do not factor through $M_2$.

\begin{teo} \label{Borromean:teo}
The isolated exceptional fillings on the Borromean rings $W$, considered up to the action of $\Iso(W)$, are listed in Table \ref{Borromean:table}.
\end{teo}

The exceptional fillings of all the 301 manifolds with $c\leq 5$ can be deduced from the theorems stated here. These 301 manifolds have either one or two cusps.

While completing this paper, we have been informed that Dunfield has recently classified and recognised the exceptional fillings of all the one-cusped hyperbolic manifolds in the much wider $c\leq 9$ census. He discovered in particular that there are almost 206,000 exceptional fillings overall \cite{Du}.

\section{Proofs}
The proofs of the Theorems \ref{main:teo},  \ref{exc:6:teo}, \ref{exc:N:teo}, \ref{exc:7:teo}, and \ref{Borromean:teo}
follow the strategy outlined in the introduction. To detect the isolated exceptional fillings we use the python script {\tt find\_exceptional\_fillings.py}, that can be downloaded from \cite{M} to be used on any multi-cusped hyperbolic manifold. The script was already used in \cite{MaPeRo}, and some of its routines have been rewritten in a more efficient way to increase its speed. We refer to \cite{MaPeRo} for a detailed tutorial to the script. 

As explained there, the script produces two finite lists of fillings: a list of ``probably isolated exceptional fillings'', and a list of ``probably hyperbolic fillings'' that typically contains only closed fillings. To conclude rigorously with the proof we need to prove that indeed all the members of the first list are not hyperbolic, and those of the second are hyperbolic (in the unlucky case where the second list contains cusped manifolds, one should run the program again on each). We have been able to accomplish that in all the steps of the proofs. As in \cite{MaPeRo}, we have used \emph{Recognizer} and \emph{Regina} for the first task, and another python script for the second: the script \emph{\tt
search\_geometric\_solution.py} is designed to determine the hyperbolicity of a long list of fillings of the same manifold $M$, using retriangulations and finite coverings. We could complete all cases using finite covers of order $\leq 11$. Few cases needed the degree 11 coverings.

Following this strategy, all the exceptional fillings of $M_i$ and $N_i$ with $i\leq 5$ can be classified directly. The program needs few seconds for $i\leq 3$, few minutes for $i=4$, and few hours for $i=5$. To attack $M_6$ we use the formidable symmetries of $M_6$ and note that, thanks to Agol and Lackenby's 6 theorem \cite{Ag, La}, an exceptional filling of $M_6$ factors either through $M_5$ or through $M_6(-2)$. So to complete the classification for $M_6$ we only needed to analyze the 5-cusped manifold $M_6(-2)$. Alternatively, we can directly classify all the fillings on $M_6$ in few days of computer time. We used a similar approach for the manifolds $N_3, \ldots, N_6$.

Classifying directly the exceptional fillings of $M_7$ takes too much computer time, so
we attacked it using the 6 theorem again: we only needed to analyse $M_7(\alpha)$ with $\alpha \in \{-3, -2, -1, -\frac 12, 0, \frac 12, 1, 2, \infty\}$. Note that, due also to the presence of fewer symmetries, we have to consider many more cases for $M_7$ than for $M_6$. The slopes $0$ and $\infty$ are exceptional. The slopes $-1$ and $1$ have already been analysed and give $M_6$ and $N_6$. We ran the program on $M_7(\alpha)$ for each remaining $\alpha$, and after a few days the computation was done.

We have also used another script to regenerate from the sometimes small amount of information found (the exceptional fillings that do not factor) all the isolated exceptional fillings of $M_6$ and $M_7$, by acting with the isometry groups of $M_6$ and $M_7$.

\setlength{\tabcolsep}{6 pt}

\begin{table}
\begin{center}
\begin{tabular}{cll}
\hline
\alza $\infty$ & $S^3$ & $\{e\}$ \\
\addlinespace[.1 cm]
\hdashline 
\addlinespace[.1 cm]
\alza $\pm 4$ & $\seifdue{D}2121
        \bigu 0110\seifdue{D}2131$ & $\matZ_4$  \\ 
\addlinespace[.1 cm]
\hdashline 
\addlinespace[.1 cm]
\alza $\pm 3$ & $\seiftre{S^2}31314{-3}$ & $\matZ_3$ \\
\addlinespace[.1 cm]
\hdashline 
\addlinespace[.1 cm]
\alza $\pm 2$ & $\seiftre{S^2}21415{-4}$ & $\matZ_2$ \\
\addlinespace[.1 cm]
\hdashline 
\addlinespace[.1 cm]
\alza $\pm 1$ & $\seiftre{S^2}21317{-6}$ & $\{e\}$ \\
\addlinespace[.1 cm]
\hdashline 
\addlinespace[.1 cm]
\alza $0$ & $T_{\tiny{\matr 31{-1}0}}$ & $\matZ$\\
\addlinespace[.2 cm]
\hline
\addlinespace[.2 cm]
\end{tabular}
\end{center}

\nota{The exceptional fillings of the figure-eight knot complement $M_1$. We also show their first homology groups.}
\label{M1_exc:table}
\end{table}

\begin{table}
\begin{center}
\begin{tabular}{cll}
\hline
\alza $\infty$ & $D \times S^1$ & $\matZ$ \\
\addlinespace[.1 cm]
\hdashline 
\addlinespace[.1 cm]
\alza $0$ & $\big(P\times S^1\big)/_{\tiny{\matr 0110}}$ & $\matZ^2$ \\
\addlinespace[.1 cm]
\hdashline 
\addlinespace[.1 cm]
\alza $1$ & $\seifdue{D}2131$ & $\matZ$ \\
\addlinespace[.1 cm]
\hdashline 
\addlinespace[.1 cm]
\alza $2$ & $\seifdue{D}2141$ & $\matZ \times \matZ_2$ \\
\addlinespace[.1 cm]
\hdashline 
\addlinespace[.1 cm]
\alza $3$ & $\seifdue{D}3131$ & $\matZ \times \matZ_3$ \\
\addlinespace[.1 cm]
\hdashline 
\addlinespace[.1 cm]
\alza $4$ & $\seifdue{D}2121
        \bigu 0110\seifuno{A}21$ & $\matZ \times \matZ_4$ \\ 
\addlinespace[.2 cm]
\hline
\addlinespace[.2 cm]
\alza $\left(\dfrac 32, 5\right)$ & $\seiftre{S^2}21313{-5}$ & $\matZ_{15}$ \\
\addlinespace[.2 cm]
\hdashline 
\addlinespace[.2 cm]
\alza $\left(\dfrac 43, 5\right)$ & 
$\seifdue{D}2121 \bigu 21{-1}{-1} \seifdue{D}2131$ & $\matZ_{20}$ \\
\addlinespace[.2 cm]
\hdashline 
\addlinespace[.2 cm]
\alza $\left(\dfrac 52, \dfrac 72 \right)$ & 
$\seifdue{D}2131 \bigu 110{-1} \seifdue{D}2131$ & $\matZ_{35}$ \\
\addlinespace[.2 cm]
\hdashline 
\addlinespace[.1 cm]
\alza $(-2,-2)$ & 
$\seifdue{D}2121 \bigu 12{0}{-1} \seifdue{D}2131$ & $\matZ_2^2$ \\
\addlinespace[.2 cm]
\hline
\addlinespace[.2 cm]
\end{tabular}
\end{center}

\nota{The isolated exceptional fillings of the Whitehead link complement $M_2$ that do not factor through $M_1$, up to the action of $\Iso(M_2)$. We also show their first homology groups.}
\label{M2_exc:table}
\end{table}

\begin{table}
\begin{center}
\begin{tabular}{cll}
\hline
\alza $\infty$ & $A \times S^1$ & $\matZ^2$ \\
\addlinespace[.1 cm]
\hdashline 
\addlinespace[.1 cm]
\alza $0$ & $\seifdue{D}2131\bigu 11{1}0 P \times S^1$ & $\matZ^2$ \\ 
\addlinespace[.1 cm]
\hdashline 
\addlinespace[.1 cm]
\alza $1$ & $\seifuno{A}21$ & $\matZ^2$ \\
\addlinespace[.1 cm]
\hdashline 
\addlinespace[.1 cm]
\alza $2$ & $\seifuno{A}31$ & $\matZ^2$ \\
\addlinespace[.1 cm]
\hdashline 
\addlinespace[.1 cm]
\alza $3$ & 
$\seifuno{A}21\bigu 0110\seifuno{A}21$ & $\matZ^2$ \\ 
\addlinespace[.2 cm]
\hline
\addlinespace[.1 cm]
\begin{tabular}{c} 
\alzamolto $\left(4, \dfrac 12\right)$ \\
\alzamolto $\left(\dfrac 32, \dfrac 52\right)$
\end{tabular}
& 
$\seifdue D2131 \bigu 11{1}0 \seifuno A21$ & $\matZ$ \\
\addlinespace[.1 cm]
\hline
\addlinespace[.2 cm]
$\left(\dfrac 83, \dfrac 32, \dfrac 32\right)$ &
    $\seifdue{D}2121
        \bigcup\nolimits_{{\tiny{\matr 21{-1}{-1}}}\phantom{\Big|}\!\!}
        \seifdue{D}2131\phantom{\Big|}$ & $\matZ_{20}$ \\ 
\addlinespace[.2 cm]
\hdashline 
\addlinespace[.2 cm]
$\left( \dfrac 52, \dfrac 53, \dfrac 53 \right)$ &
    $\seifdue{D}2121
        \bigcup\nolimits_{{\tiny{\matr {3}{1}{-2}{-1}}}\phantom{\Big|}\!\!}
        \seifdue{D}2131\phantom{\Big|}$ & $\matZ_{16}$ \\ 
\addlinespace[.2 cm]
\hdashline 
\addlinespace[.1 cm]
$\left( -2, -2, -2 \right)$ &
    $\seifdue{D}2121
        \bigcup\nolimits_{{\tiny{\matr 23{-1}{-2}}}\phantom{\Big|}\!\!}
        \seifdue{D}2131\phantom{\Big|}$ & $\matZ_4$ \\
\addlinespace[.2 cm]
\hdashline
\addlinespace[.2 cm]
\alza $\left(4,\dfrac 32, \dfrac 32\right)$ 
& $T_{{\tiny{\matr {-3}1{-1}0}}\phantom{\Big|}}\phantom{\Big|}$ & $\matZ \times \matZ_5$ \\ 
\addlinespace[.2 cm]
\hdashline 
\addlinespace[.2 cm]
$\left(5, 5, \dfrac 12\right)$ &
    $\seifuno{A}21\big/_{{\tiny{\matr 0110}}\phantom{\Big|}}\phantom{\Big|}$ & $\matZ \times \matZ_3$ \\
\addlinespace[.2 cm]
\hdashline 
\addlinespace[.1 cm]
\begin{tabular}{c}
\alzamolto $\left(4, 4, \dfrac 23\right)$ \\
\alzamolto $\left( \dfrac 73, \dfrac 73, \dfrac 32 \right)$ 
\end{tabular}
&
    $\seifuno{A}21\big/_{{\tiny{\matr 1110}}\phantom{\Big|}}\phantom{\Big|}$ & $\matZ \times \matZ_5$\\
\addlinespace[.1 cm]
\hdashline 
\addlinespace[.2 cm]
$\left(\dfrac 52, \dfrac 52, \dfrac 43\right)$ &
    $\seifuno{A}21\big/_{{\tiny{\matr 2110}}\phantom{\Big|}}\phantom{\Big|}$ & $\matZ \times \matZ_7$\\
\addlinespace[.2 cm]
\hline
\addlinespace[.2 cm]
\end{tabular}
\end{center}

\nota{The isolated exceptional fillings of the magic manifold $M_3$ that do not factor through $M_2$, up to the action of $\Iso(M_3)$. We also show their first homology groups.}
\label{M3_exc:table}
\end{table}

\begin{table}
\begin{center}
\begin{tabular}{cll}
\hline
 $\infty$ & \alza $P\times S^1$ & $\matZ^3$ \\
\addlinespace[.1 cm]
\hdashline 
\addlinespace[.1 cm]
   $0$ & \alza $\big(P\times S^1\big) \bigu 0110 \seifuno A21$ & $\matZ^3$ \\
\addlinespace[.2 cm]
\hline
\addlinespace[.2 cm]
\alza $\left(-2, -2, -2, -2\right)$
& $\seifdue{D}212{1} \bigu {3}4{-2}{-3} \seifdue{D}2131$ & $\matZ_2^2$ \\
\addlinespace[.2 cm]
\hline
\addlinespace[.2 cm]
\end{tabular}
\end{center}

\nota{The isolated exceptional fillings of $M_4$ that do not factor through $M_3$, up to the action of $\Iso(M_4)$. We also show their first homology groups.}
\label{M4_exc:table}
\end{table}

\begin{table}
\begin{center}
\begin{tabular}{cll}
\hline
 $\infty$ & \alza $\big(P\times S^1\big) \bigu 0110 \big(P\times S^1\big)$ & $\matZ^4$ \\
\addlinespace[.2 cm]
\hline
\addlinespace[.2 cm]
\alza $\left(-2, -2, -2, -2, -2\right)$
& $\seifdue{D}212{1} \bigu {4}5{-3}{-4} \seifdue{D}2131$ & $\matZ_4$ \\
\addlinespace[.2 cm]
\hdashline 
\addlinespace[.2 cm]
\alza $\left(-2, -2, \dfrac 13, 3, \dfrac 13 \right)$
& $\seifuno A 2{1} \bigb 120{-1}$ & $\matZ \times \matZ_6$ \\
\addlinespace[.2 cm]
\hline
\addlinespace[.2 cm]
\end{tabular}
\end{center}

\nota{The isolated exceptional fillings of $M_5$ that do not factor through $M_4$, up to the action of $\Iso(M_5)$. We also show their first homology groups.}
\label{M5_exc:table}
\end{table}

\begin{table}
\begin{center}
\begin{tabular}{cll}
\hline
 $\infty$ & \alza $\big(P\times S^1\big) \bigu 0110 \big(P\times S^1\big) \bigu 0110 \big(P\times S^1\big)$ & $\matZ^5$ \\
\addlinespace[.2 cm]
\hline
\addlinespace[.2 cm]
\alza $(-2,-2,\boldsymbol\cdot,2,2)$ & $\seifdue D313{2} \bigu 1110 \big(P\times S^1\big)$ & $\matZ^2 \times \matZ_3$ \\
\addlinespace[.2 cm]
\hdashline 
\addlinespace[.2 cm]
\alza $\left(-2,-\dfrac{1}{2},\boldsymbol\cdot,\dfrac{1}{2},2\right)$ & 
$\big(P\times S^1\big) 
\bigcup\nolimits_{{\tiny{\matr 0110}}}^{{\tiny{\matr 0110}}} \big(P\times S^1\big)$ 
& $\matZ^3$
\\
\addlinespace[.2 cm]
\hline
\addlinespace[.2 cm]
\alza $\left(2,2,\dfrac 12, -3,-2\right)$ & 
$\seifdue D2152$
& $\matZ$
\\ 
\addlinespace[.2 cm]
\hdashline 
\addlinespace[.2 cm]
\alza $\left(2,2,\dfrac 12, -4,-2\right)$ & 
$\seifdue D213{1} \bigu 1110 \seifuno A21$
& $\matZ$
\\
\addlinespace[.2 cm]
\hdashline 
\addlinespace[.2 cm]
\alza $\left(2,2,-\dfrac{1}{2},-2,-\dfrac{3}{2}\right)$ & 
$\seifdue D214{1} \bigu 1110 \seifuno A21$
& $\matZ \times \matZ_2$
\\
\addlinespace[.2 cm]
\hdashline 
\addlinespace[.2 cm]
\alza $\left(\dfrac 12, -2,-\dfrac{1}{2},-2,\dfrac{1}{2}\right)$ & 
$\big(P\times S^1\big) \bigb 1211$
& $\matZ^2 \times \matZ_2$
\\
\addlinespace[.2 cm]
\hline
\addlinespace[.2 cm]
\end{tabular}
\end{center}

\nota{The non-closed isolated exceptional fillings of $M_6$ that do not factor through $M_5$, up to the action of $\Iso(M_6)$. The dot indicates a cusp that is not filled: it appears when the fillings are not along consecutive components on the link. We also show the first homology groups.}
\label{M6_exc:table}
\end{table}

\begin{table}
\begin{center}
\begin{tabular}{cll}
\hline
\addlinespace[.2 cm]
\begin{tabular}{c}
\alzamolto $\left(-3,-\dfrac{2}{3},-2,\dfrac{1}{2},2,\dfrac{1}{2}\right)$ \\
\alzamolto $\left(-3,-\dfrac{1}{2},-3,\dfrac{1}{2},2,\dfrac{1}{2}\right)$ 
\end{tabular}
& 
$\seiftre{S^2} 21318{-7}$
& $\matZ_2$
\\
\addlinespace[.2 cm]
\hdashline 
\addlinespace[.2 cm]
\begin{tabular}{c}
\alzamolto $\left(-3,-2,-\dfrac{2}{3},-2,2,\dfrac{1}{2}\right)$ \\
\alzamolto $\left(-3,-2,-\dfrac{1}{2},-2,2,\dfrac{1}{3}\right)$ 
\end{tabular}
& 
$\seiftre {S^2}2131{11}{-9}$ 
&
$\{e\}$ 
\\
\addlinespace[.2 cm]
\hdashline 
\addlinespace[.2 cm]
\begin{tabular}{c}
\alzamolto $\left(-4,-\dfrac{2}{3},-2,\dfrac{1}{2},2,\dfrac{1}{2}\right)$ \\
\alzamolto $\left(-4,-\dfrac{1}{2},-3,\dfrac{1}{2},2,\dfrac{1}{2}\right)$ \\
\alzamolto $\left(-\dfrac{4}{3},-\dfrac{1}{2},2,\dfrac{1}{2},2,-\dfrac{1}{3}\right)$ 
\end{tabular} & 
$\seiftre {S^2}21416{-5}$
&
$\matZ_2^2$
\\
\addlinespace[.1 cm]
\hdashline 
\addlinespace[.2 cm]
\begin{tabular}{c}
\alzamolto $\left(-3,-\dfrac{2}{3},-2,\dfrac{1}{3},2,\dfrac{1}{2}\right)$ \\ 
\alzamolto $\left(-3,-\dfrac{1}{2},-3,\dfrac{1}{2},2,\dfrac{1}{3}\right)$ 
\end{tabular}
&
$\seiftre{S^2} 21515{-4}$
&
$\matZ_5$
\\ 
\addlinespace[.2 cm]
\hdashline 
\addlinespace[.2 cm]
\begin{tabular}{c}
\alzamolto $\left(-5,-\dfrac{2}{3},-2,\dfrac{1}{2},2,\dfrac{1}{2}\right)$ \\
\alzamolto $\left(-5,-\dfrac{1}{2},-3,\dfrac{1}{2},2,\dfrac{1}{2}\right)$ \\
\alzamolto $\left(-\dfrac{5}{4},-\dfrac{1}{2},2,\dfrac{1}{2},2,-\dfrac{1}{3}\right)$ 
\end{tabular}
&  
$\seiftre {S^2} 313{1}5{-4}$ 
& $\matZ_6$ \\ 
\addlinespace[.1 cm]
\hdashline 
\addlinespace[.1 cm]
\begin{tabular}{c}
\alzamolto $\left(-4,-\dfrac{2}{3},-2,\dfrac{1}{3},2,\dfrac{1}{2}\right)$ \\
\alzamolto $\left(-4,-\dfrac{1}{2},-3,\dfrac{1}{3},2,\dfrac{1}{2}\right)$ 
\end{tabular}
& 
$\seiftre {S^2}31414{-3}$ 
& $\matZ_8$
\\
\addlinespace[.2 cm]
\hline
\addlinespace[.2 cm]
\end{tabular}
\end{center}

\nota{The closed isolated exceptional fillings of $M_6$ that do not factor through $M_5$, up to the action of $\Iso(M_6)$. We also show their first homology groups. (Part I.)}
\label{M6b_exc:table}
\end{table}

\begin{table}
\begin{center}
\begin{tabular}{cll}
\hline
\begin{tabular}{c}
\alzamolto $\left(-6,-\dfrac{2}{3},-2,\dfrac{1}{2},2,\dfrac{1}{2}\right)$  \\
\alzamolto $\left(-6,-\dfrac{1}{2},-3,\dfrac{1}{2},2,\dfrac{1}{2}\right)$ \\
\alzamolto $\left(-\dfrac{6}{5},-\dfrac{1}{2},2,\dfrac{1}{2},2,-\dfrac{1}{3}\right)$ 
\end{tabular}
&
    $\seifdue{D}2121
        \bigu 0110
        \seifdue{D}2141$ 
&
$\matZ_2 \times \matZ_4$
\\
\addlinespace[.1 cm]
\hdashline 
\addlinespace[.1 cm]
\begin{tabular}{c}
\alzamolto $\left(-5,-\dfrac{2}{3},-2,\dfrac{1}{3},2,\dfrac{1}{2}\right)$ \\
\alzamolto $\left(-5,-\dfrac{1}{2},-3,\dfrac{1}{3},2,\dfrac{1}{2}\right)$ 
\end{tabular}
& 
    $\seifdue{D}2131
        \bigu 0110
        \seifdue{D}2131$ 
&
$\matZ_{11}$ \\ 
\addlinespace[.1 cm]
\hdashline 
\addlinespace[.1 cm]
\begin{tabular}{c}
\alzamolto $\left(-\dfrac{5}{2},-2,-\dfrac{2}{5},2,2,\dfrac{1}{2}\right)$ \\
\alzamolto $\left(-\dfrac{3}{2},-2,-\dfrac{2}{3},2,2,\dfrac{1}{2}\right)$ 
\end{tabular}
& 
    $\seifdue{D}2131
        \bigu 11{0}{-1}
        \seifdue{D}2132$ 
&
$\matZ_{37}$ \\
\addlinespace[.1 cm]
\hdashline 
\addlinespace[.2 cm]
\alza $\left(-\dfrac{5}{3},-2,-\dfrac{1}{2},2,2,\dfrac{1}{3}\right)$ & 
    $\seifdue{D}2131
        \bigu 01{1}{-1}
        \seifdue{D}213{2}$ 
&
$\matZ_{31}$ 
\\  
\addlinespace[.1 cm]
\hdashline 
\addlinespace[.1 cm]
\begin{tabular}{c}
\alzamolto $\left(-4,-\dfrac{2}{3},-2,\dfrac{1}{4},2,\dfrac{1}{2}\right)$ \\
\alzamolto $\left(-4,-\dfrac{1}{2},-3,\dfrac{1}{4},2,\dfrac{1}{2}\right)$  
\end{tabular}
&
    $\seifdue{D}2121
        \bigu 01{1}0 
        \seifdue{D}3131$ 
&
$\matZ_2 \times \matZ_6$
\\  
\addlinespace[.2 cm]
\hdashline 
\addlinespace[.2 cm]
$\left(-3,-\dfrac{1}{2},2,2,2,-\dfrac{1}{2}\right)$ & 
    $\seifdue{D}2121
        \bigu 21{-1}{-1}
        \seifdue{D}3131$   
&
$\matZ_{24}$ \\
\addlinespace[.2 cm]
\hdashline
\addlinespace[.2 cm]
\begin{tabular}{c}
\alzamolto $\left(-4,-\dfrac{3}{4},-2,\dfrac{1}{2},2,\dfrac{1}{2}\right)$ \\
\alzamolto $\left(-4,-\dfrac{1}{2},-4,\dfrac{1}{2},2,\dfrac{1}{2}\right)$ 
\end{tabular}
 & 
    $\seifdue{D}2121
        \bigu {1}20{-1}
        \seifdue{D}214{1}\phantom{\Big|}$ 
&
$\matZ_2^3$ \\  
\addlinespace[.1 cm]
\hdashline 
\addlinespace[.2 cm]
$\left(-3,-\dfrac{1}{2},-3,\dfrac{1}{3},2,\dfrac{1}{3}\right)$ & 
    $\seifdue{D}2121
        \bigu {1}20{-1} 
        \seifdue{D}3131$ & $\matZ_{12}$ \\  
\addlinespace[.2 cm]
\hline 
\addlinespace[.2 cm]
\end{tabular}
\end{center}

\nota{The closed isolated exceptional fillings of $M_6$ that do not factor through $M_5$, up to the action of $\Iso(M_6)$. We also show their first homology groups. (Part II.)}
\label{M6c_exc:table}
\end{table}

\begin{table}
\begin{center}
\begin{tabular}{cll}
\hline
\addlinespace[.2 cm]
\begin{tabular}{c}
\alzamolto $\left(-4,-\dfrac{1}{2},-3,\dfrac{1}{2},2,\dfrac{1}{3}\right)$ \\
\alzamolto $\left(-\dfrac{4}{3},-\dfrac{1}{2},3,\dfrac{1}{2},2,-\dfrac{1}{3}\right)$ 
\end{tabular}
& 
    $\seifdue{D}2131
        \bigu {1}20{-1} 
        \seifdue{D}2131$ 
& $\matZ_{10}$ \\  
\addlinespace[.1 cm]
\hdashline 
\addlinespace[.2 cm]
$\left(-3,-\dfrac{2}{3},-3,\dfrac{1}{2},2,\dfrac{1}{2}\right)$ & 
    $\seifdue{D}2121
        \bigu {2}3{-1}{-2}
        \seifdue{D}214{1}\phantom{\Big|}$ & $\matZ_2 \times \matZ_4$ \\  
\addlinespace[.1 cm]
\hdashline 
\addlinespace[.2 cm]
$\left(-3,-\dfrac{2}{3},-2,\dfrac{1}{2},2,\dfrac{1}{3}\right)$ & 
    $\seifdue{D}2131
        \bigu {2}3{-1}{-2} 
        \seifdue{D}2131$ & $\matZ_9$ \\  
\addlinespace[.1 cm]
\hdashline 
\addlinespace[.2 cm]
$\left(-\dfrac{3}{2},-\dfrac{1}{2},2,\dfrac{2}{3},2,-\dfrac{1}{2}\right)$ & 
    $\seifdue{D}2131
        \bigu {3}4{-2}{-3}
        \seifdue{D}2131$ & $\matZ_8$ \\
\addlinespace[.2 cm]
\hdashline
\addlinespace[.2 cm]
$(-2,-2,-2,-2,-2,-2)$ & 
    $\seifdue{D}2121
        \bigu {5}6{-4}{-5}
        \seifdue{D}2131$ & $\matZ_2^2$ \\
\addlinespace[.1 cm]
\hdashline 
\addlinespace[.1 cm]
\begin{tabular}{c}
\alzamolto $\left(-4,-2,-\dfrac{2}{3},-2,2,\dfrac{1}{2}\right)$ \\ 
\alzamolto $\left(-3,-3,-\dfrac{1}{2},-3,2,\dfrac{1}{2}\right)$ 
\end{tabular}
& 
$\seifuno{A}21\big/_{{\tiny{\matr {-1}110}}\phantom{\Big|}}\phantom{\Big|}$ 
& $\matZ$ \\
\addlinespace[.1 cm]
\hdashline 
\addlinespace[.2 cm]
$\left(-3,-3,\dfrac{1}{2},2,2,\dfrac{1}{2}\right)$ & 
$\seifuno{A}21\big/_{{\tiny{\matr {1}30{-1}}}\phantom{\Big|}}\phantom{\Big|}$ 
& $\matZ \times \matZ_5$ \\ 
\addlinespace[.2 cm]
\hdashline 
\addlinespace[.2 cm]
\begin{tabular}{c}
\alzamolto $\left(-3,-2,-\dfrac{1}{3},3,2,\dfrac{1}{3}\right)$ \\
\alzamolto $\left(-\dfrac{3}{2},-2,\dfrac{1}{2},\dfrac{3}{2},2,-\dfrac{1}{2}\right)$ 
\end{tabular} 
&
$M_1 \bigu {-1}0{1}{1} \seifdue D2{1}2{1}$ & $\matZ_2^2$ \\
\addlinespace[.2 cm]
\hdashline 
\addlinespace[.2 cm]
\alza $\left(-3,-\dfrac{3}{2},-\dfrac{1}{2},2,2,\dfrac{1}{3}\right)$ & 
$M_1 \bigu {-1}110 \seifdue D2121$ & $\matZ_4$ \\
\addlinespace[.2 cm]
\hline
\addlinespace[.2 cm]

\end{tabular}
\end{center}

\nota{The closed isolated exceptional fillings of $M_6$ that do not factor through $M_5$, up to the action of $\Iso(M_6)$. We also show their first homology groups. (Part III.)}
\label{M6d_exc:table}
\end{table}

\begin{table}
\begin{center}
\begin{tabular}{cll}
\hline
\addlinespace[.1 cm]
\alza $\infty$ & $A \times S^1$ & $\matZ^2$ \\
\addlinespace[.1 cm]
\hdashline 
\addlinespace[.1 cm]
\alza $0$ & $M_1 \bigu 0110 \big(P\times S^1\big) $ & $\matZ^2$ \\
\addlinespace[.2 cm]
\hline 
\addlinespace[.1 cm]
\alza $(2,2,2)$ & $\seifuno A 2{1} \bigb 2311$ & $\matZ \times \matZ_3$\\
\addlinespace[.2 cm]
\hline 
\addlinespace[.1 cm]
\end{tabular}
\end{center}
\nota{The isolated exceptional fillings of $N_3$ that do not factor, up to the action of $\Iso(N_3)$. We also show their first homology groups.}
\label{N3:table}
\end{table}

\begin{table}
\begin{center}
\begin{tabular}{cll}
\hline
\addlinespace[.1 cm]
\alza $\infty$ & $P\times S^1$ & $\matZ^3$ \\
\addlinespace[.1 cm]
\hdashline 
\addlinespace[.1 cm]
\alza $0$ & $M_2 \bigu 0110 \big(P\times S^1\big) $ & $\matZ^3$ \\
\addlinespace[.2 cm]
\hline
\addlinespace[.2 cm]
\alza $(2,2,2,2)$ & $\seifuno A 2{1} \bigb 3411$ & $\matZ \times \matZ_4$ \\
\addlinespace[.2 cm]
\hline 
\addlinespace[.1 cm]
\end{tabular}
\end{center}
\nota{The isolated exceptional fillings of $N_4$ that do not factor, up to the action of $\Iso(N_4)$. We also show their first homology groups.}
\label{N4:table}
\end{table}

\begin{table}
\begin{center}
\begin{tabular}{cll}
\hline
\addlinespace[.1 cm]
\alza $\infty$ & $\big(P\times S^1\big) \bigu 0110 \big(P\times S^1 \big)$ & $\matZ^4$ \\
\addlinespace[.1 cm]
\hdashline 
\addlinespace[.1 cm]
\alza $0$ & $M_3 \bigu 0110 \big(P\times S^1\big) $ & $\matZ^4$ \\
\addlinespace[.2 cm]
\hline
\addlinespace[.2 cm]
\alza $\left(-4, -\dfrac 32, -\dfrac 32, -4, -\dfrac 12\right)$ & $M_1 \bigu {-1}21{-1} \seifdue D2121$ & $\matZ_2 \times \matZ_2$ \\
\addlinespace[.2 cm]
\hdashline
\addlinespace[.2 cm]
\alza $\left( -\dfrac 32, -\dfrac 32, -\dfrac 32, -\dfrac 32, -\dfrac 32 \right)$
& $\seifdue D2{1}31 \bigu {4}{5}{-3}{-4} \seifdue D2131$ & $\matZ_7$ \\
\addlinespace[.2 cm]
\hdashline
\addlinespace[.2 cm]
\alza $(2,2,2,2,2)$ & $\seifuno A 2{1} \bigb 4511$ & $\matZ \times \matZ_5$ \\
\addlinespace[.2 cm]
\hline 
\addlinespace[.1 cm]
\end{tabular}
\end{center}
\nota{The isolated exceptional fillings of $N_5$ that do not factor, up to the action of $\Iso(N_5)$. We also show their first homology groups.}
\label{N5:table}
\end{table}

\begin{table}
\begin{center}
\begin{tabular}{cll}
\hline
\alza $\infty$ & \alza $\big(P\times S^1\big) \bigu 0110 \big(P\times S^1\big) \bigu 0110 \big(P\times S^1\big)$ & $\matZ^5$ \\
\addlinespace[.1 cm]
\hdashline 
\addlinespace[.1 cm]
\alza $0$ & $M_4 \bigu 0110 \big(P\times S^1\big) $ & $\matZ^5$ \\
\addlinespace[.2 cm]
\hline
\addlinespace[.2 cm]
\alza $\left(-\dfrac 12, -4, -\dfrac 12, -4, -\dfrac 12\right)$ & 
$M_2 \bigu {-1}21{-1} \seifdue D2121$ & $\matZ \times \matZ_2^2$ \\
\addlinespace[.2 cm]
\hdashline
\addlinespace[.2 cm]
\alza $\left(-\dfrac 12, -3, -\dfrac 23, -3, -\dfrac 12\right)$ & 
$ M_2 \bigu {-1}110 \seifdue D2121$  & $\matZ \times \matZ_4$ \\
\addlinespace[.2 cm]
\hline
\addlinespace[.2 cm]
\alza $\left(-5, -\dfrac 12, -3, -\dfrac 12, -3, -\dfrac 12\right)$ & 
$\seifdue D2131 \bigu 0110 \seifdue D2131$
& $\matZ_{11}$ \\
\addlinespace[.2 cm]
\hdashline 
\addlinespace[.2 cm]
\alza $\left(-4, -\dfrac 12, -3, -\dfrac 12, -3, -\dfrac 12\right)$ & 
$\seiftre {S^2} 31414{-3}$
& $\matZ_8$ \\
\addlinespace[.2 cm]
\hdashline 
\addlinespace[.2 cm]
\alza $\left(-3, -3, -\dfrac 12, -3, -3, -\dfrac 12\right)$ &
$\seifuno A21 \bigb 140{-1}$
& $\matZ \times \matZ_4$ \\
\addlinespace[.2 cm]
\hdashline 
\addlinespace[.2 cm]
\alza $\left(-3, -\dfrac 12, -3, -\dfrac 12, -3, -\dfrac 12\right)$ &
$\seiftre {S^2} 21515{-4}$ & $\matZ_5$ \\
\addlinespace[.2 cm]
\hdashline 
\addlinespace[.2 cm]
\alza $\left(-3, -\dfrac 12, 2, -\dfrac 12, 2, -\dfrac 12 \right)$ &
$\seifuno A 2{1} \bigb 1211$ & $\matZ \times \matZ_2$ \\
\addlinespace[.2 cm]
\hdashline 
\addlinespace[.2 cm]
\alza $(2,2,2,2,2,2)$ &
$\seifuno A 2{1} \bigb 5611$ & $\matZ \times \matZ_6$ \\
\addlinespace[.2 cm]
\hline 
\addlinespace[.1 cm]
\end{tabular}
\end{center}

\nota{The isolated exceptional fillings of $N_6$ that do not factor, up to the action of $\Iso(N_6)$. We also show their first homology groups.}
\label{N6:table}
\end{table}

\begin{table}
\begin{center}
\begin{tabular}{cll}
\hline
\addlinespace[.2 cm]
\alza $\infty$ & $(D \times S^1) \# (D \times S^1)$ & $\matZ^2$ \\
\addlinespace[.2 cm]
\hdashline
\addlinespace[.2 cm]
\alza $0$ & 
$\big(P\times S^1\big) 
\bigcup\nolimits_{{\tiny{\matr 0{-1}{-1}0}}}^{{\tiny{\matr 0110}}} \big(P\times S^1\big)$ 
& $\matZ^3$ \\
\addlinespace[.2 cm]
\hline
\addlinespace[.2 cm]
$(-2,-2)$ & 
$\seifdue D2121 \bigu 1211 \seifuno A21$
& $\matZ \times \matZ_2^2$ \\
\addlinespace[.2 cm]
\hline
\addlinespace[.2 cm]
$\left(-\dfrac 43, -3, -2 \right)$ & 
$\seifdue D2121 \bigu 21{-1}{-1} \seifdue D2141$
& $\matZ_2 \times \matZ_{12}$ \\
\addlinespace[.2 cm]
\hdashline
\addlinespace[.2 cm]
$\left(-3, -\dfrac 32, -2 \right)$ & 
$\seiftre {S^2}31313{-4}$
& $\matZ_3 \times \matZ_6$ \\
\addlinespace[.2 cm]
\hline
\addlinespace[.2 cm]
\end{tabular}
\end{center}
\nota{The isolated exceptional fillings of the Borromean rings complement $W$ that do not factor through $M_1$, up to the action of $\Iso(W)$. We also show their first homology groups.}
\label{Borromean:table}
\end{table}

\begin{table}
\begin{center}
\begin{tabular}{cll}
\hline
\addlinespace[.2 cm]
\alza $\infty$ & \alza \begin{tabular}{c} $\big(P\times S^1\big) \bigu 0110 \big(P\times S^1\big) \bigu 0110$ \\ $\big(P\times S^1\big) \bigu 0110 \big(P\times S^1\big)$ \end{tabular} & $\matZ^6$ \\
\addlinespace[.2 cm]
\hdashline
\addlinespace[.1 cm]
\alza $0$ & $M_5 \bigu 0110 \big(P\times S^1\big) $ & $\matZ^6$ \\
\addlinespace[.2 cm]
\hline
\addlinespace[.1 cm]
\begin{tabular}{c}
\alzamolto
$\left(\dfrac 12, 2, -\dfrac 12, -\dfrac 32, -2, -\dfrac 12\right)$ \\
\alzamolto
$\left(-\dfrac 12, -2, \dfrac 12, \dfrac 12, -2, -\dfrac 12\right)$
\end{tabular}
& 
$M_2 \bigu {-1}0{1}{1} \seifdue D 2121 $
& $\matZ \times \matZ_2^2$
\\
\addlinespace[.1 cm]
\hline
\addlinespace[.1 cm]
\end{tabular}
\end{center}
\nota{The non-closed isolated exceptional fillings of $M_7$ that do not factor, up to the action of $\Iso(M_7)$. We also show their first homology groups.}
\label{M7a:table}
\end{table}

\begin{table}
\begin{center}
\begin{tabular}{cll}
\hline
\addlinespace[.2 cm]
$\left(-3, -2, -\dfrac 12, 2, 2, 2, -\dfrac 12\right)$ & 
$\seifdue D 2121 \bigu 21{-1}{-1} \seifdue D 3141$
& $\matZ_{28}$
\\
\addlinespace[.1 cm]
\hdashline 
\addlinespace[.1 cm]
$\left(-3, -2, \dfrac 12, 2, 2, \dfrac 12, -2\right)$
& $\seifdue D 2121 \bigu 21{-1}{-1} \seifdue D2151$ 
& $\matZ_{16}$ \\
\addlinespace[.1 cm]
\hdashline 
\addlinespace[.2 cm]
$\left(-2, -\dfrac 12, 2, -\dfrac 12, -\dfrac 12, 2, -\dfrac 12 \right)$ & 
$\seifdue D2121 \bigu 34{-2}{-3} \seifdue D2141 $
& $\matZ_2^3$ \\
\addlinespace[.2 cm]
\hdashline 
\addlinespace[.2 cm]
$\left( -\dfrac 12, -\dfrac 12, -\dfrac 12, -\dfrac 12, -\dfrac 12, -\dfrac 12, -\dfrac 12 \right) $ & 
$\seifdue D 2131 \bigu 67{-5}{-6} \seifdue D2131$ & $\matZ_5$ \\
\addlinespace[.1 cm]
\hdashline 
\addlinespace[.1 cm]
\begin{tabular}{c}
\alzamolto $\left(-\dfrac 32, -2, -\dfrac 12, 2, 2, -\dfrac 12, -2 \right)$ \\
\alzamolto $\left(-\dfrac 32, -2, -\dfrac 12, 2, 2, \dfrac 12, 2 \right)$
\end{tabular}
& $\seifuno A 2{1} \bigb 1211$ & $\matZ \times \matZ_2$\\
\addlinespace[.1 cm]
\hdashline 
\addlinespace[.2 cm]
$(2,2,2,2,2,2,2) $ & 
$\seifuno A21 \bigb 6711$
& $\matZ \times \matZ_7$ \\
\addlinespace[.2 cm]
\hdashline 
\addlinespace[.1 cm]
$\left(-3, -2, -3, \dfrac 12, 2, 2, \dfrac 12 \right)$ & $\seifuno A 3{1} \bigb 130{-1}$ & 
$\matZ \times \matZ_9$ \\
\addlinespace[.1 cm]
\hdashline 
\addlinespace[.1 cm]
\begin{tabular}{c}
\alzamolto $\left(-\dfrac 32, -\dfrac 32, -2, \dfrac 12, 2, \dfrac 12, -2 \right)$ \\
\alzamolto $\left( -\dfrac 32, -2, -\dfrac 12, 2, \dfrac 32, \dfrac 12, -2 \right)$
\end{tabular} 
& $\big(M_2(-2)\big) \bigu {-1}0{1}{1} \seifdue D2121$ & $\matZ_2^3$ \\
\addlinespace[.1 cm]
\hdashline 
\addlinespace[.2 cm]
$\left(-3, -\dfrac 32, -2, -\dfrac 12, 2, 2, \dfrac 13 \right)$ & $\big(M_2(-2)\big) \bigu {-1}{1}{1}0 \seifdue D2121$ & 
$\matZ_2 \times \matZ_4$ \\
\addlinespace[.2 cm]
\hline 
\addlinespace[.1 cm]
\end{tabular}
\end{center}
\nota{The closed isolated exceptional fillings of $M_7$ that do not factor, up to the action of $\Iso(M_7)$. We also show their first homology groups.}
\label{M7b:table}
\end{table}

\setlength{\tabcolsep}{6 pt}


\begin{thebibliography}{99}
\bibitem{M} \bysame, personal web page, {\tt people.dm.unipi.it/martelli/research.html}
\end{thebibliography}
\end{document}